\documentclass[journal,9pt,twoside,web]{ieeecolor}
\usepackage{generic}

\usepackage{ifpdf}
\ifpdf
\fi

\usepackage[T1]{fontenc}
\usepackage[utf8]{inputenc} % or \usepackage[utf8x]{inputenc} for more characters
\usepackage{graphicx,color} % for pdf, bitmapped graphics files
\usepackage{epsfig} % for postscript graphics files
\usepackage{times} % assumes new font selection scheme installed
\usepackage{amsmath} % assumes amsmath package installed
\usepackage{amssymb}  % assumes amsmath package installed
\usepackage{caption}
\usepackage{subcaption}
\usepackage{stackrel}
\usepackage{eucal}
\usepackage{url}
\usepackage{multirow}
\usepackage{cite}
\usepackage{algorithm}
\usepackage[noend]{algpseudocode}
\usepackage{pgf}
\usepackage{pgfplots}
\usepgfplotslibrary{fillbetween}
\usepackage{tikz}
\usetikzlibrary{arrows,automata}
\usetikzlibrary{intersections}
\usetikzlibrary{calc}
\usetikzlibrary{shapes}
\usetikzlibrary{positioning}
\usetikzlibrary{decorations.text}
\pgfdeclarelayer{background}
\pgfdeclarelayer{foreground}
\pgfsetlayers{background,main,foreground}
\usetikzlibrary{calc,patterns,decorations.pathmorphing,decorations.markings}

\definecolor{mycol}{RGB}{19,48,128}
\usepackage{hyperref} 
\hypersetup{
	colorlinks = true,
	citecolor = mycol,
	linkcolor = mycol,
	urlcolor = black
}
\usepackage{cleveref}

\usepackage{tcolorbox}
\tcbuselibrary{skins}
\usetikzlibrary{calc,patterns,decorations.pathmorphing,decorations.markings}
\usetikzlibrary{shadows}
\usetikzlibrary{shadows.blur}
\usepackage{xcolor}

\usepackage[nosfdefault]{comicneue}
\newtcolorbox{stickynotebox}{
  enhanced,
  sharp corners,
  colback=yellow!25,
  frame hidden,
  boxrule=0pt,
  left=6pt,
  right=6pt,
  top=6pt,
  bottom=6pt,
  drop shadow={opacity=0.12, shadow xshift=2.2mm, shadow yshift=-2.2mm}
}

\graphicspath{{figures/}{biography/}}

\def\zb{\boldsymbol{z}}
\def\ab{\boldsymbol{a}}
\def\db{\boldsymbol{d}}

\def\yb{\boldsymbol{y}}

\def\R{\mathbb{R}}
\def\Wb{\boldsymbol{W}}
\def\xb{\boldsymbol{x}}
\def\pb{\boldsymbol{p}}
\def\bb{\boldsymbol{b}}
\def\sb{\boldsymbol{s}}
\def\vb{\boldsymbol{v}}
\def\ub{\boldsymbol{u}}
\def\Mb{\boldsymbol{M}}
\def\Qb{\boldsymbol{Q}}
\def\bla{\boldsymbol{\lambda}}
\def\gb{\boldsymbol g}
\def\Fb{\boldsymbol{F}}
\def\Ab{\boldsymbol{A}}
\def\Cb{\boldsymbol{C}}
\def\Pb{\boldsymbol{P}}
\newtheorem{lem}{Lemma}
\newtheorem{definition}{Definition}
\newtheorem{theorem}{Theorem}
\newtheorem{assumption}{Assumption}
\newtheorem{proposition}{Proposition}
\newtheorem{corollary}{Corollary}
\newtheorem{example}{Example}
\newtheorem{remark}{Remark}
\crefname{lem}{Lemma}{Lemmas}
\Crefname{lem}{Lemma}{Lemmas}
\crefname{definition}{Definition}{Definitions}
\Crefname{definition}{Definition}{Definitions}
\crefname{theorem}{Theorem}{Theorems}
\Crefname{theorem}{Theorem}{Theorems}
\crefname{assumption}{Assumption}{Assumptions}
\Crefname{assumption}{Assumption}{Assumptions}
\crefname{proposition}{Proposition}{Propositions}
\Crefname{proposition}{Proposition}{Propositions}
\crefname{corollary}{Corollary}{Corollaries}
\Crefname{corollary}{Corollary}{Corollaries}
\crefname{example}{Example}{Examples}
\Crefname{example}{Example}{Examples}
\crefname{remark}{Remark}{Remarks}
\Crefname{remark}{Remark}{Remarks}
\crefname{section}{Section}{Sections}
\Crefname{section}{Section}{Sections}
\crefname{appendix}{Appendix}{Appendices}
\Crefname{appendix}{Appendix}{Appendices}
\newcommand{\VI}{\operatorname{VI}}

\newif\ifshownew \shownewtrue
\newif\ifshowold \showoldtrue
\newif\ifcolorchanges \colorchangestrue
\newcommand{\new}[1]{\ifshownew{\ifcolorchanges\color{red}\fi\relax#1}\fi}
\newcommand{\old}[1]{\ifshowold{\ifcolorchanges\color{gray}\fi\relax#1}\fi}

\begin{document}

	%\title{\LARGE Approximate Multiparametric Generalized Nash Equilibrium Problems and Data-driven Multi-agent Explicit Model Predictive Control}
	\title{\LARGE A Hierarchical Variational Inequality Approach to Generalized Nash Equilibrium Seeking with Nonlinear Shared Constraints}

    \author{Tatiana Tatarenko, Puya Latafat, and Alberto Bemporad %\IEEEmembership{Fellow, IEEE}
    \thanks{T. Tatarenko is with the Technical University of Darmstadt, Darmstadt, Germany (e-mail: \texttt{\scriptsize tatiana.tatarenko@tu-darmstadt.de}).
    P.~Latafat and A. Bemporad are with the IMT School for Advanced Studies Lucca, Lucca, Italy (e-mail: \texttt{\scriptsize puya.latafat@imtlucca.it}; \texttt{\scriptsize alberto.bemporad@imtlucca.it}).}
    \thanks{The work was funded by the Deutsche Forschungsgemeinschaft
(DFG, German Research Foundation). Project number: 528033031.
    This work was also funded by the European Union (ERC Advanced Research Grant COMPACT, No. 101141351). Views and opinions expressed are however those of the authors only and do not necessarily reflect those of the European Union or the European Research Council. Neither the European Union nor the granting authority can be held responsible for them.
    }}
	
\maketitle
\thispagestyle{empty}
	
%%%%%%%%%%%%%%%%%%%%%%%%%%%%%%%%%%%%%%%%%%%%%%%%%%%%%%%%%%%%%%%%%%%%%%%%%%%%%%%%
\begin{abstract}
We study the computation of variational generalized Nash equilibria (v-GNE) in games with strongly monotone pseudo-gradients and nonlinear convex shared constraints. 
We show that the v-GNE problem is exactly equivalent to a hierarchical variational inequality, whose lower level is given by the gradient of a smooth convex penalty of the constraint violation over the local action sets. Building on this reformulation, we propose a hierarchical regularization scheme that operates entirely in the primal space and requires only projections onto the local action sets. In the centralized setting, we prove last-iterate squared distance convergence to the v-GNE with rate $O(1/t^{\scriptscriptstyle 2/3})$ under a local error bound condition, which holds, e.g., under Slater's condition, improved to $O(1/t)$ under a global error bound condition. We then extend the scheme to the partial-decision information setting, where the agents communicate over a network and maintain local estimates of the joint action. We derive a distributed algorithm with a constant stepsize operating in the primal space and converging to the v-GNE with rate $O(1/t)$. 
\end{abstract}

\begin{keywords}
 Variational generalized Nash equilibrium, hierarchical variational inequality, distributed algorithms, nonlinear shared constraints
\end{keywords}

\section{Introduction}
%\PL{\textit{I am marking edit with \texttt{\textbackslash new} (red) and \texttt{\textbackslash old} (gray). Both are shown by default. 
%To only see the new text, put \texttt{\textbackslash showoldfalse} in the preamble.}}
Generalized Nash equilibrium problems (GNEPs) model noncooperative interactions in which self-interested agents are coupled not only through their cost functions but also through shared constraints on their joint actions~\cite{fischer2014generalized}. Such problems arise in a variety of engineering applications, including electricity markets and demand-side management~\cite{elmark,chen2014autonomous}, charging coordination of electric vehicles~\cite{ma2011decentralized}, traffic and transportation networks~\cite{dafermos1980traffic}, and communication systems~\cite{scutari2012monotone}. In many of these settings, the shared constraints represent physical limits of a common infrastructure, such as line ratings, transformer loadings, or voltage-stability margins in power networks, and are therefore inherently \emph{nonlinear}~\cite{low2014convex}. Moreover, the agents are typically not coordinated by a central operator and possess only local information: each agent knows its own cost and local constraints and can exchange information only with its neighbors over a communication network~\cite{molzahn2017survey}.

Among the generalized Nash equilibria, a particularly relevant subclass is formed by the \emph{variational} GNE (v-GNE), which correspond to solutions of a variational inequality (VI) defined by the game pseudo-gradient over the coupled feasible set~\cite{facchinei2007generalized}. The v-GNE are economically meaningful, as all agents face the same shadow price for the shared resources~\cite{KULKARNI2012}, and, under strong monotonicity of the pseudo-gradient, a v-GNE exists and is unique~\cite{FaccPang1}. The VI characterization suggests solving the problem by projection-type methods, such as projected gradient play or the extragradient method~\cite{Korp76,BakPolyak74,NesterovScrimali,Lan2022}. However, these methods require projections onto the coupled feasible set. When the shared constraints are nonlinear, such projections are generally not available in closed form and require solving an auxiliary optimization problem at every iteration. Furthermore, since the coupled feasible set depends on the joint action profile, these projections cannot be computed locally by the agents.

For this reason, the majority of existing GNE seeking algorithms rely on the Karush--Kuhn--Tucker (KKT) conditions of the v-GNE and operate in the lifted primal--dual space, in which the shared constraints are handled through Lagrange multipliers~\cite{tat_kam_TAC,ran2024distributed}. In the partial-decision information setting, where each agent maintains estimates of the actions of the others and communicates only with its neighbors~\cite{TatNedShi21}, distributed primal--dual schemes based on operator splitting have been developed in~\cite{pavel2020distributed,bianchi2022fast}. These methods, however, require the agents to maintain, update, and reach consensus on local copies of the dual variables, which increases both the computational and the communication burden. Moreover, their convergence analysis concerns linear coupling constraints.  For nonlinear shared constraints, the lifted operator is not globally Lipschitz continuous, even if the original game's pseudo-gradient possesses this property, and, thus, constant-stepsize primal--dual methods lose their guarantees. In this case, additional mechanisms, such as adaptive or line-search stepsize rules~\cite{malitsky2020golden,MalitskyMishchenko2020}, are required. However, these rules are still applied to the lifted primal--dual operator and cannot be implemented in a distributed manner, as the line search requires the stepsize to be coordinated centrally among all agents. %Linear convergence under nonlinear coupled inequality constraints has been established in~\cite{meng2023linear}, yet in a centralized full-information setting.

In this paper, we take a different route and avoid dual variables altogether. We observe that the set of feasible joint actions coincides with the set of minimizers of a smooth convex penalty function defined for the constraint violation over the local action sets. As a consequence, the v-GNE problem is \emph{exactly} equivalent to a \emph{hierarchical} VI, in which the upper level is given by the game pseudo-gradient and the lower-level VI is given by the gradient of the penalty and involves only projections onto the local action sets. Hierarchical VIs have been studied by means of vanishing (Tikhonov-type) regularization~\cite{Kalashnikova, MSt25} (see also~\cite{Attouch1996} for minimization problems) and extragradient-type schemes~\cite{EG_HVI_Samadi}. However, all these methods are {centralized} procedures: they assume that the full operators of both levels are available to a single computing unit and require full information on the joint decision vector at every iteration. These works address general monotone lower-level operators and typically provide guarantees for ergodic averages or in terms of gap functions. In contrast, the lower-level VI in our reformulation corresponds to a convex optimization problem, and we exploit this structure to obtain last-iterate convergence rates to the v-GNE. To the best of our knowledge, the hierarchical VI reformulation has not been used before to design distributed GNE seeking algorithms that operate entirely in the primal space.

The main contributions of this paper are summarized as follows.
\begin{itemize}
    \item We propose an equivalent reformulation of the v-GNE problem with nonlinear convex shared constraints as a hierarchical VI, which does not require the introduction of dual variables.
    \item We develop a centralized hierarchical regularization scheme that performs a single projected step onto the local action sets per iteration. Under a local error bound condition on the constraint set, which holds, e.g., under Slater's condition, we prove last-iterate squared distance convergence to the v-GNE with rate $O(1/t^{\scriptscriptstyle 2/3})$. Under a global error bound condition, which holds, e.g., for linear shared constraints and polyhedral local sets, we improve this rate to $O(1/t)$.
%\PL{The rates $O(1/t^{2/3})$ and $O(1/t)$ in this bullet, in the abstract, and in the conclusion are for the squared distance $\|\ab(t)-\ab^\star\|^2$. The distance itself decays as $O(1/t^{1/3})$ and $O(1/t^{1/2})$. We should say ``squared distance'' explicitly.}
    \item We extend the scheme to the partial-decision information setting, in which the agents communicate over a connected undirected graph and maintain local estimates of the joint action. The resulting distributed algorithm uses a constant stepsize and requires no communication of dual variables. We prove its convergence to the v-GNE with rate $O(1/t)$ { under the local error bound condition (e.g., Slater's condition) and bounded action sets}.
    \item We illustrate the performance of the proposed algorithms on randomly generated instances %inspired by the operation of multi-energy systems with nonlinear network constraints 
    of multi-agent systems with nonlinear network constraints
    and on distributed problems with linear shared constraints, and compare them with the corresponding KKT-based lifted primal--dual methods~\cite{malitsky2020golden,Lan2022,bianchi2022fast}.
\end{itemize}

The remainder of the paper is organized as follows. Section~\ref{sec:problem} formulates the problem and discusses the limitations of the projection-based and primal--dual approaches. Section~\ref{sec:HVI} presents the hierarchical VI reformulation. The centralized and distributed hierarchical regularization schemes and their convergence analysis are presented in Sections~\ref{sec:centralized} and~\ref{sec:distributed}, respectively. Section~\ref{sec:simulations} provides the simulation results, and the proofs of auxiliary lemmas are collected in the Appendix.

\subsection{Notation}
We write $[N]=\{1,\ldots,N\}$ and $\R^m_+$ for the nonnegative orthant of $\R^m$. For vectors $\ub^1,\ldots,\ub^N$ of possibly different dimensions, $[\ub^1,\ldots,\ub^N]$ denotes the stacked column vector; a superscript $i$ refers to the block associated with player $i$. For $\ub,\vb\in\R^n$, $\langle\ub,\vb\rangle$ is the Euclidean inner product and $\|\ub\|$ is the Euclidean norm; for a matrix, $\|\cdot\|$ is the induced spectral norm. We denote by $\boldsymbol 0$ the zero vector, by $\boldsymbol 1_N$ the $N$-dimensional vector of ones, by $I_N$ the $N\times N$ identity matrix, by $A^\top$ the transpose of a matrix $A$, and by $\otimes$ the Kronecker product. For a symmetric matrix $A$, $A\succeq0$ means that $A$ is positive semidefinite, and $\lambda_{\min}(A)$ denotes its smallest eigenvalue. For $t\in\R$, $[t]_+=\max\{t,0\}$; for vectors, $[\cdot]_+$ is applied componentwise. For a differentiable function $\phi:\R^D\to\R$, $\nabla\phi(\ab)$ is its gradient and $\nabla_i\phi(\ab)$ is its partial gradient with respect to $\ab^i$.

For a nonempty closed convex set $X\subseteq\R^n$, $\Pi_X(\ub)=\operatorname{argmin}_{\yb\in X}\|\yb-\ub\|$ is the Euclidean projection of $\ub$ onto $X$, $\operatorname{dist}(\ub,X)=\|\ub-\Pi_X(\ub)\|$ is the distance from $\ub$ to $X$, $\EuScript N_X(\ub)=\{\vb:\langle\vb,\yb-\ub\rangle\le0\ \forall\yb\in X\}$ is the normal cone to $X$ at $\ub\in X$, and $\operatorname{int}X$ and $\operatorname{ri}X$ are the interior and the relative interior of $X$. Given $\Fb:X\to\R^n$, $\VI(X,\Fb)$ is the problem of finding $\ub\in X$ such that $\langle\Fb(\ub),\yb-\ub\rangle\ge0$ for all $\yb\in X$, and $\operatorname{SOL}(X,\Fb)$ is its solution set.

A mapping $\Fb:X\to\R^n$ is \emph{$L$-Lipschitz continuous} on $X$ if $\|\Fb(\ub)-\Fb(\vb)\|\le L\|\ub-\vb\|$, \emph{monotone} if $\langle\Fb(\ub)-\Fb(\vb),\ub-\vb\rangle\ge0$, \emph{$\mu$-strongly monotone} if $\langle\Fb(\ub)-\Fb(\vb),\ub-\vb\rangle\ge\mu\|\ub-\vb\|^2$, and \emph{$\frac1L$-cocoercive} if $\langle\Fb(\ub)-\Fb(\vb),\ub-\vb\rangle\ge\frac1L\|\Fb(\ub)-\Fb(\vb)\|^2$, for all $\ub,\vb\in X$. If $\Fb=\nabla\phi$ for a convex differentiable function $\phi$ and $\Fb$ is $L$-Lipschitz continuous on an open convex set $U$, then $\Fb$ is $\frac1L$-cocoercive on $U$~\cite[Thm.~3.1]{PerezArosVilches2021}. For nonnegative sequences, $u_t=O(v_t)$ means that there exists $c>0$ such that $u_t\le c\,v_t$ for all sufficiently large $t$.

\section{Problem Formulation}\label{sec:problem}
We consider a non-cooperative game $\Gamma$ between $N$ players. Let player $i$, $i \in [N]$, choose its local action from $\R^{d_i}$ and $J^i: \R^D\to\R$, with $D=\sum_{i=1}^Nd_i$, denote its cost function. We assume the existence of a local constraint set $\Ab^i\subseteq \R^{d_i}$ for each player $i\in[N]$ as well as of a \emph{shared coupling constraint set} $\Cb$ defined by the following system of  inequalities: 
\begin{align}\label{eq:coupled_cs}
\Cb = \{\ab\in\R^D: \gb(\ab)\le \boldsymbol{0}\},
\end{align}
where $\gb(\ab) = [g_1(\ab),\ldots,g_m(\ab)]: \R^D\to\R^m$,  $\ab = [\ab^1,\ldots,\ab^N]\in\R^D$, with $\ab^i = [a^{i,1},\ldots,a^{i,d_i}]\in\R^{d_i}$ being the action of  player $i$ and $\ab^{-i}$ denoting the joint action of all players except player $i$. %Thus, $C = \{\ab\in\R^D: \gb(\ab) = K\ab-\lb\le \boldsymbol{0}\}$ with $K=\sum_{i=1}^N \boldsymbol{e}_i\otimes K^i$ and $\lb = [\lb^1,\ldots,\lb^N]\in\R^{n}$, $n=\sum_{i=1}^Nn_i$.
Thus, the feasible action set in game $\Gamma$ is $\Ab\cap \Cb$, where $\Ab = \Ab^1\times\ldots\times\Ab^N$, and the coupled action set of agent $i$ is $\Ab^i\cap{\Cb}^i(\ab^{-i})$, where ${\Cb}^i(\ab^{-i}) = \{\ab^i\in\R^{d_i}: \gb(\ab)\le \boldsymbol{0}\}$. Throughout this work, we assume that the individual action sets $\Ab^i$, $i \in [N]$, possess a structure that permits efficient computation of Euclidean projections.

A \emph{generalized Nash equilibrium} (GNE) in  game $\Gamma$ represents a joint action from which no player has any incentive to unilaterally deviate.
\begin{definition}\label{def:GNE}
 A joint action $\ab^\star\in  \Ab\cap \Cb$ in the game $\Gamma=\Gamma(N, \{\Ab^i\}, \{J^i\}, \Cb)$ with the coupling constraints $\Cb$ is called a \emph{generalized Nash equilibrium} (GNE) if $J^i(\ab^{\star,i},\ab^{\star,-i})\le J^i(\ab^i,\ab^{\star,-i})$ for any $i\in[N]$ and $\ab^i\in\Ab^i\cap \Cb^i(\ab^{\star,-i})$.
 \end{definition}

%We are interested in designing a \emph{payoff-based algorithm}, which converges to a generalized Nash equilibrium in $\Gamma$. The payoff-based information structure implies that  agent $i$, $i \in [N]$, does not know the functional form of $J_i$ or $\gb$ but for a played action $\ab$, it can observe their values. 
We consider convex games as defined below. 
\begin{assumption}\label{assum:convex}
 The game under consideration is \emph{convex}. Namely,
{\begin{enumerate}
\item[(i)] for all $i\in[N]$ the cost function $J^i(\ab^i, \ab^{-i})$ is defined on $\R^{D}$, continuously differentiable in $\ab$ and convex in $\ab^i$ for  fixed $\ab^{-i}$;
\item[(ii)] the set $\Ab$ is closed and convex,  ${g_j}:\R^D\to\R$ is convex {and continuously differentiable} for any ${j\in[m]}$, and the set ${\Ab\cap \Cb}$ is nonempty.
\end{enumerate}} %The coupling constraint function $\gb:\R^{Nd}\to\R^n$ is continuously differentiable and has convex coordinates $g_i(\ab), i\in[n]$.
\end{assumption}
%\begin{assumption}\label{assum:Slater}
%  There exists $\bar\ab \in \Ab$ such that
%\(
% g_i(\bar\ab) < 0, \qquad i = 1,\dots,m .
%\)
%\end{assumption}
Given the differentiability of the cost functions, we can define the game's pseudo-gradient mapping below.
\begin{definition}\label{def:mapping}
The mapping $\Mb = [\Mb^1, \ldots, \Mb^N]:\R^{D}\to\R^{D}$, referred to as  the \emph{pseudo-gradient} of  $\Gamma(N, \{\Ab^i\}, \{J^i\}, \Cb)$, is defined by
 \begin{align}\label{eq:gamemapping}
 &\Mb^i(\ab) = [M^{i,1}(\ab), \ldots, M^{{i},d_i}(\ab)], \cr
 &M^{i,j}(\ab)= \frac{\partial J^i(\ab)}{\partial a^{i,j}}, \; \ab\in  {\R^D}, j\in[d_i], i\in[N].
 \end{align}
 \end{definition}
 %\PL{I also removed the transpose. If I understood the notation, a stacked vector is denoted using \(w = [w^1,\ldots, w^N]\), so there is no need for transpose for \(M^i\).}

In general, characterizing and computing generalized Nash equilibria (GNE) of the game $\Gamma$ is challenging \cite{facchinei2007generalized}. In this work, we focus on the set of \emph{variational generalized Nash equilibria} (v-GNE), which can be characterized as the solution to the variational inequality defined below.
\begin{definition}\label{def:VGNE}
    The set of \emph{variational generalized Nash equilibria  (v-GNE)} in the game $\Gamma$ is defined as the set of solutions to the following variational inequality problem:
    \begin{align}\label{eq:GVI}
    \nonumber
    \mbox{Find } &\ab^\star\in  \Ab\cap \Cb: \, \langle\Mb(\ab^\star),\ab-\ab^\star\rangle\ge 0, \\
    &\mbox{ for all }\ab\in\Ab\cap \Cb .%\tag{$VI{(\Ab\cap C, \Mb)}$}
\end{align}
\end{definition}
The next lemma states the relation between GNE and v-GNE (see \cite{facchinei2007generalized} for the proof). 
\begin{lem}\label{lem:GNEandVGNE}
    Under Assumption~\ref{assum:convex}, the set of generalized Nash equilibria in  game $\Gamma$ contains the set of variational generalized Nash equilibria, i.e., v-GNE $\subseteq$ GNE\footnote{Without loss of generality we assume that for any set $\EuScript S$ the following holds: $\emptyset\subseteq \EuScript S$, where $\emptyset$ denotes the empty set.}.   
\end{lem}

{Convexity of each $J^i$ in its own action $\ab^i$, as required in Assumption~\ref{assum:convex}{(i)}, does not imply monotonicity of the pseudo-gradient $\Mb$. Monotonicity is an additional property of the joint mapping and is central to the analysis of variational inequality problems.}
In this work we are interested in games with a specific structure which allows for guaranteeing existence and uniqueness of a v-GNE, and are amenable to efficient distributed solution schemes.
{For this reason we make the following additional strong monotonicity assumption on $\Mb$.}
\begin{assumption}\label{assum:strMon}
  The pseudo-gradient $\Mb$ is $L_{\Mb}$-Lipschitz continuous and $\mu$-strongly monotone over $\Ab$\footnote{Note that Assumption~\ref{assum:strMon} implies that $J^i(\ab^{i},\ab^{-i})$ is {$\mu$-strongly convex in $\ab^i\in\Ab^i$ for every fixed $\ab^{-i}\in\prod_{j\neq i}\Ab^j$, which strengthens the convexity in $\ab^i$} (the first part of Assumption~\ref{assum:convex}).}.   
\end{assumption}
{Strong convexity of each $J^i$ in its own action $\ab^i$ alone does not imply strong monotonicity of $\Mb$, since the latter also depends on the coupling between players. A similar strong monotonicity and Lipschitz continuity assumption is imposed in~\cite{bianchi2022fast}, against which we compare in \cref{sec:simulations}.}
%\PL{Should we cite them and others here that assume strong monotonicity?}
Taking into account that $\Ab\cap \Cb\subseteq \Ab$ and the result on existence of solutions to VI (see, for example,~\cite{facchinei2007generalized}), we conclude that under Assumptions~\ref{assum:convex}--\ref{assum:strMon} the game $\Gamma$ possesses a unique v-GNE, i.e., a solution to~\eqref{eq:GVI}.
\begin{remark}\label{rem:VI}
A straightforward approach for computing a v-GNE is to apply gradient-play or projected-gradient-type methods to the variational inequality~\eqref{eq:GVI}. Such methods are particularly efficient when the feasible set $\Ab\cap \Cb$ admits a projection operator that can be computed in closed form or at low computational cost. However, when the coupling constraints defining $\Cb$ are nonlinear, evaluating the projection onto $\Ab\cap \Cb$ is, in general, nontrivial and may itself require solving a separate optimization problem at each iteration.

This difficulty becomes more pronounced in the generalized Nash setting, where the constraint set is coupled through the actions of all players. As a result, the projection onto $\Ab\cap \Cb$ cannot be carried out locally by each agent using only its own variables; instead, it depends on the joint action profile, which is not directly available under distributed information structures.
Consequently, it is unclear how to efficiently implement projection-based VI methods in a decentralized manner when agents have access only to partial information about other players' decisions. This motivates the development of alternative approaches that avoid explicit projections onto the coupled constraint set and rely instead on more tractable, distributed-friendly mechanisms.
\end{remark}

{To address such alternative approaches, in the following analysis we will require the distance to the feasible set $\Ab\cap \Cb$ to be controlled by the constraint violation, at least on bounded sets.}
\begin{assumption}\label{assum:local-eb}
{The feasible set $\Ab\cap \Cb$ admits a local error bound on $\Ab$, i.e., for every bounded set $B\subseteq\R^D$ there exists $\kappa_B>0$ such that
$\operatorname{dist}(\ab,\Ab\cap \Cb)\le\kappa_B\,\|[\gb(\ab)]_+\|$, $\forall\ab\in B\cap\Ab$.
}
\end{assumption}
{\Cref{assum:local-eb} holds, in particular, under Slater's condition, i.e., if there exists $\bar\ab\in\Ab$ such that $g_j(\bar\ab)<0$ for all $j\in[m]$ (see \cref{lem:slater-eb} in \cref{app:slater}). It also holds, even globally on $\Ab$ and without a Slater point, when $\gb$ is affine and $\Ab$ is polyhedral, by Hoffman's error bound (see \cref{rem:AffineConst} below). Note that \cref{assum:local-eb} is not implied by \cref{assum:convex} alone: for $\Ab=\R$ and $g(a)=a^2$, we have $\Ab\cap \Cb=\{0\}$ and $\operatorname{dist}(a,\Ab\cap \Cb)=|a|$, which cannot be bounded by $\kappa_B[g(a)]_+=\kappa_B a^2$ near $a=0$.}

\begin{remark}\label{rem:lVI}
{A natural alternative to projecting onto $\Ab\cap \Cb$ is to work with the KKT system of~\eqref{eq:GVI}, which requires a constraint qualification. Slater's condition}, together with the convexity of the game in Assumption~\ref{assum:convex}, allows us to formulate the following equivalent Karush--Kuhn--Tucker (KKT) conditions characterizing a v-GNE $\ab^\star$: there exists $\bla^\star\in\R^m$ (a common Lagrange multiplier vector) such that
\begin{align}\label{eq:KKT}
\boldsymbol 0 &\in \Mb(\ab^\star) + \nabla \gb(\ab^\star)^\top \bla^\star
+ \EuScript N_{\Ab}(\ab^\star), \cr
\ab^\star &\in \Ab,\qquad
\gb(\ab^\star) \le \boldsymbol 0, \cr
\bla^\star &\ge \boldsymbol 0, \cr
\langle \bla^\star,\gb(\ab^\star) \rangle &= 0,
\end{align}
The KKT conditions~\eqref{eq:KKT} are equivalent to the following variational inequality on the lifted (primal--dual) space:
\begin{align}
\text{Find } (\ab^\star,\bla^\star)\in \Ab\times\R_+^m \text{ such that }  \, \forall (\ab,\bla)\in \Ab\times\R_+^m\cr\quad
\left\langle
\begin{bmatrix}
\Mb(\ab^\star)+\nabla\gb(\ab^\star)^\top\bla^\star\\
-\gb(\ab^\star)
\end{bmatrix},
\begin{bmatrix}
\ab-\ab^\star\\
\bla-\bla^\star
\end{bmatrix}
\right\rangle
\ge 0.
\label{eq:liftedVI}
\end{align}

It is worth noting that in the case of \emph{linear coupling constraints} $\gb(\ab)\le \boldsymbol 0$, {the KKT conditions~\eqref{eq:KKT} characterize the v-GNE without Slater's condition, provided that $\Ab$ is polyhedral or $\operatorname{ri}\Ab\cap \Cb\neq\emptyset$ \cite[Cor. 23.8.1]{rockafellar1970convex}}. Moreover, in this case the operator
\[
F(\ab,\bla)=
\begin{bmatrix}
\Mb(\ab)+\nabla\gb(\ab)^\top\bla\\
-\gb(\ab)
\end{bmatrix}
\]
is monotone and Lipschitz continuous under Assumption~\ref{assum:strMon}. {However, strong monotonicity of $\Mb$ does not provide strong monotonicity of the lifted operator.} {Nevertheless, mere monotonicity and Lipschitz continuity} enable the use of efficient algorithms and their distributed variants for solving the monotone VI~\eqref{eq:liftedVI}; see, e.g.,~\cite{TatKamECC25,bianchi2022fast}. However, such methods require working in the lifted space and performing updates of the dual variable $\bla$, which may be computationally and communication-wise expensive in distributed implementations. 

Furthermore, for general nonlinear coupling constraints, the mapping $F$ is not globally Lipschitz continuous on $\Ab\times\R_+^m$. This limitation necessitates additional central coordination or adaptive stepsize strategies, such as line-search schemes, to ensure convergence of such methods~\cite{malitsky2020golden}.

In contrast to lifted-space methods that require explicit introduction and updates of the dual variables $\bla$ and rely on solving a variational inequality with a potentially non-Lipschitz continuous operator, in this work we develop a primal approach. Our method operates directly in the space of the primal actions $\ab$ and does not require maintaining or updating Lagrange multipliers, thereby avoiding the computational and communication overhead associated with dual variables and lifted formulations.
\end{remark}

\section{Hierarchical VI Reformulation}\label{sec:HVI}
Motivated by the observations made in the previous section (see Remarks~\ref{rem:VI},~\ref{rem:lVI}), we now derive an equivalent reformulation of the considered GNE problem that is expressed entirely in terms of the primal decision variables and thus avoids the introduction of dual variables and the associated lifted-space variational inequality formulation.

Let us define any \emph{penalty profile} $\psi$ from the following class.
\begin{definition}\label{def:profile}
A function $\psi:[0,\infty)\to[0,\infty)$ is an \emph{admissible penalty profile} if $\psi$ is continuously differentiable{ with $\psi'(0)=0$}, $\psi$ is convex and nondecreasing, $\psi(s)=0 \iff s=0$, and there exist $S_0>0$ and $c_\psi>0$ such that
$\psi(s) \ge c_\psi\, s\,\min(s,S_0)$, $\forall s\ge0$,
i.e., $\psi$ grows at least quadratically near the origin, $\psi(s)\ge c_\psi s^2$ for $s\le S_0$, and at least linearly beyond the threshold $S_0$, $\psi(s)\ge c_\psi S_0\, s$ for $s> S_0$.
\end{definition}
%\PL{I think we also need $\psi'(0)=0$. Without it, $\psi(s)=s$ is admissible with $c_\psi=1/S_0$, and then $\Psi=\sum_j[g_j]_+$ need not be differentiable at points where some $g_j$ vanishes.}

\begin{example}[Quadratic profile]\label{ex:quad}
The profile $\psi(s)=\frac12 s^2$ satisfies Definition~\ref{def:profile} with $c_\psi=\frac12$ for \emph{every} choice of $S_0>0$.
\end{example}

\begin{example}[Huber profile]\label{ex:huber}
For $\omega>0$, the capped profile
\[
\psi_\omega(s) = \begin{cases} \frac12 s^2, & 0\le s\le \omega,\\ \omega s-\frac12\omega^2, & s>\omega,\end{cases}
\]
coincides with the quadratic profile on $[0,\omega]$ and grows only linearly beyond the cap $\omega$. %; capping the profile this way keeps $\Pb=\nabla\Psi$ from inheriting an arbitrarily large Lipschitz constant when $\gb$ is steep. It satisfies Definition~\ref{def:profile} with $c_\psi=\frac12$ and $S_0=M$: for $s\le M$, $\psi_M(s)=\frac12s^2=\frac12s\min(s,M)$ (equality); for $s>M$, $\psi_M(s)-\frac12Ms=Ms-\frac12M^2-\frac12Ms=\frac12M(s-M)\ge0$. Unlike the quadratic profile of Example~\ref{ex:quad}, $\psi_M$ satisfies no purely quadratic bound $\psi_M(s)\ge cs^2$ valid for \emph{all} $s\ge0$ with a fixed $c>0$ -- which is exactly why Definition~\ref{def:profile} allows a linear regime beyond $S_0$.
\end{example}

Given an admissible penalty profile $\psi$, let us define the penalty function for the coupling constraints as follows:
\begin{align}\label{eq:pen}
    \Psi(\ab) {:=} \sum_{j=1}^m \psi([g_j(\ab)]_+),
\end{align}
together with its gradient
\begin{align}\label{eq:LL}
\Pb(\ab) {:=} \nabla \Psi(\ab), \quad \Pb^i(\ab) {:=} \nabla_i\Psi(\ab).
\end{align}
{The smoothness properties of $\Pb$ used in the following analysis are addressed in \cref{app:penalty},  \cref{lem:P-Lip}.}

\begin{proposition}\label{prop:pen-basic}
{Let Assumption~\ref{assum:convex}(ii) hold.} For any admissible penalty profile $\psi$: (i) $\Psi(\ab)=0 \iff \ab\in \Cb$; (ii) $\Psi$ is convex, so that $\Pb=\nabla\Psi$ is monotone over $\Ab$.
\end{proposition}
\begin{proof}
(i) By Definition~\ref{def:profile}, $\psi(s)=0\iff s=0$, so $\Psi(\ab)=\sum_j\psi([g_j(\ab)]_+)=0$ iff $[g_j(\ab)]_+=0$ for all $j\in[m]$, iff $\gb(\ab)\le\boldsymbol0$, iff $\ab\in \Cb$.
(ii) Since $g_j$ is convex over $\Ab$ for every $j\in[m]$ (Assumption~\ref{assum:convex}{(ii)}), $[g_j]_+$ is convex; since $\psi$ is convex and nondecreasing, the composition $\psi([g_j(\cdot)]_+)$ is convex, and a sum of convex functions is convex, so $\Psi$ is convex. Convexity of $\Psi$ implies monotonicity of $\Pb=\nabla\Psi$ over $\Ab$.
\end{proof}

{The} set $\Ab\cap \Cb$ is {equal} to the solution set $\operatorname{SOL}(\Ab,\Pb)$ {of} the following ${\VI}(\Ab,\Pb)$:
 \[
    \text{find } \ab \in \Ab \text{ such that }
    \langle \Pb(\ab), \yb - \ab \rangle \geq 0
    \quad \forall \yb \in \Ab.
    \]
{Indeed, since $\Psi$ is convex and continuously differentiable, as shown in \cref{lem:P-Lip}, $\operatorname{SOL}(\Ab,\Pb)=\operatorname{argmin}_{\Ab}\Psi$. By \cref{prop:pen-basic}, since $\Psi$ is nonnegative and vanishes exactly on $\Cb$ and $\Ab\cap \Cb\neq\emptyset$ as in \cref{assum:convex}(ii), this set equals $\Ab\cap \Cb$.}
%\PL{I defined \texttt{\textbackslash VI} as \texttt{\textbackslash operatorname\{VI\}} in the preamble and used it only in ``$\VI(\Ab,\Pb)$'' above, for comparison with $\VI(\Ab,\Pb)$. It sets VI upright like $\operatorname{SOL}$. The paper currently mixes \texttt{\$\VI(} and \texttt{\textbackslash operatorname\{VI\}}. If you like it, we can use \(\VI\) throughout.}
Thus, the original v-GNE problem is equivalent to the following \emph{hierarchical VI}:
    \begin{align}\label{eq:HVI}
    \text{find } &\ab^\star \in S_L \text{ such that }
    \langle \Mb(\ab^\star), \yb - \ab^\star \rangle \geq 0\cr
    &\forall \yb \in S_L= \operatorname{SOL}(\Ab,\Pb).
    \end{align}
    {Since $S_L=\Ab\cap \Cb$, problem~\eqref{eq:HVI} is exactly $\VI(\Ab\cap \Cb,\Mb)$. Representing its feasible set as $\operatorname{SOL}(\Ab,\Pb)$ allows us to construct methods using evaluations of $\Pb$ and projections onto $\Ab$, without projecting onto $\Ab\cap \Cb$.}
%\PL{I added the short argument for $\operatorname{SOL}(\Ab,\Pb)=\Ab\cap C$ and rephrased the last sentence. I kept the exact equivalence explicit and highlighted its algorithmic benefit that the representation allows us to construct methods using projections onto $\Ab$ only. Does this make the motivation clearer?}

We make the following assumptions regarding the constraint set and, thus, the mapping $\Pb$. 
\begin{assumption}\label{assum:Lip}
    $\Pb$ is $L_{\Pb}$-Lipschitz continuous over { an open convex set $U\subseteq\R^D$ containing $\Ab$, with $L_{\Pb}>0$}.
    %\PL{See my comment after \eqref{eq:regTh1}.}
\end{assumption}
\begin{remark}
Assumption~\ref{assum:Lip} is standard in the literature on hierarchical and bilevel variational inequalities. In particular, convergence and complexity analyses of algorithms for hierarchical VI problems typically require the operator defining the lower-level variational inequality to be monotone and Lipschitz continuous; see, e.g., \cite{EG_HVI_Samadi} and the references therein. In our setting, the operator $\Pb$ arises as the gradient of the smooth penalty function $\Psi$, and Assumption~\ref{assum:Lip} can be verified directly from the regularity properties of the constraint mapping $\gb$, as discussed below.
{If $\psi'$ is Lipschitz continuous, as for the profiles of \cref{ex:quad,ex:huber}, \cref{lem:P-Lip} with $X=U$ yields \cref{assum:Lip} in each of the following cases.
    \begin{itemize}
        \item The set $\Ab$ is compact and $\gb$ is twice continuously differentiable. Choose a bounded open convex neighborhood $U$ of $\Ab$; compactness of $\overline U$ gives the bounds required by \cref{lem:P-Lip}.
        \item $\nabla\gb$ is bounded and Lipschitz continuous over an open convex neighborhood $U$ of $\Ab$, and $\psi'$ is bounded, as for the Huber profile, for which the constant $\bar\psi$ of \cref{lem:P-Lip} equals $\sqrt m\,\omega$.
        \item $\gb$ is affine, so that the constant $L_{\gb}$ of \cref{lem:P-Lip} equals $0$, and we may take $U=\R^D$.
    \end{itemize}
The constant $L_{\Pb}$ used in the stepsize restrictions below is the Lipschitz bound on $U$.}
%\PL{I replaced the two conditions by \cref{lem:P-Lip} and stated the bounds on a neighborhood of $\Ab$ to match \cref{assum:Lip}. The boundedness of $\psi'$ in the second case controls the term involving the variation of $\nabla\gb$.}

However, in contrast to the general hierarchical VI setting considered in the existing literature, the lower-level problem in our reformulation possesses a specific structure. In particular, since $\Pb=\nabla\Psi$, where $\Psi$ is a convex differentiable penalty function, the lower-level variational inequality $\VI(\Ab,\Pb)$ is equivalent to a convex optimization problem. As will be shown later, this additional structure can be exploited in the convergence analysis of the proposed algorithmic scheme. In particular, unlike the extragradient-based methods developed for general hierarchical VIs~\cite{EG_HVI_Samadi}, which typically provide convergence guarantees for ergodic averages, the proposed approach allows us to establish explicit convergence rate guarantees for the last iterate.    
\end{remark}

\section{Centralized Hierarchical Regularization Scheme}\label{sec:centralized}

Before turning to the distributed setting, we first solve the hierarchical VI~\eqref{eq:HVI} by a centralized scheme, i.e., one in which the mappings $\Pb$ and $\Mb$ are evaluated at the true, commonly available joint action $\ab(t)$. This scheme both serves as a performance benchmark for, and forms the basis of the derivation of, the distributed algorithm of Section~\ref{sec:distributed}.

Let each player $i$ start from some action $\ab^i(0)\in\Ab^i$ and update it using first-order information, i.e., the local gradients $\Mb^i$, together with the lower-level operator $\Pb^i$, as summarized in Algorithm~\ref{alg:central}.

\begin{algorithm}[t]
\caption{Centralized hierarchical regularization scheme}
\label{alg:central}
\begin{algorithmic}[1]
\State \textbf{Initialization:} Each agent chooses $\ab^i(0)\in\Ab^i$. The stepsizes $\gamma,\eta>0$ and the regularization sequence $\{\varepsilon_t\}\subseteq(0,\infty)$ are common to all agents.
\For{$t=0,1,2,\ldots$}
\State Each agent $i\in[N]$ updates its action as
\begin{align*}
   \ab^i(t+1)
=
\Pi_{\Ab^i}\Bigl[
\ab^i(t) - \gamma \Pb^i(\ab(t)) - \eta\varepsilon_t \Mb^i(\ab(t))
\Bigr].
\end{align*}
\EndFor
\end{algorithmic}
\end{algorithm}

Thus, the sequence of joint actions $\{\ab(t)\}$ is generated by
\begin{equation}\label{eq:algC}
\ab(t+1)
=
\Pi_{\Ab}\Bigl[
\ab(t) - \gamma \Pb(\ab(t)) - \eta\varepsilon_t \Mb(\ab(t))
\Bigr],
\end{equation}
which can be regarded as a single gradient step towards solving the regularized $\VI(\Ab,\Pb+\tau_t\Mb)$, with effective regularization parameter $\tau_t=\frac{\eta}{\gamma}\varepsilon_t$. Unlike the distributed scheme developed in Section~\ref{sec:distributed}, here every agent has access to the full joint action $\ab(t)$, so that no consensus or estimation mechanism is required; Algorithm~\ref{alg:central} may thus be viewed as the idealized limit of the distributed Algorithm~\ref{alg:1} introduced below, in which consensus is instantaneous, i.e., every agent's local estimate of the joint action coincides with $\ab(t)$ at every $t$.

We analyze convergence of~\eqref{eq:algC} using two complementary techniques. The first, presented next, follows the classical vanishing-regularization (Tikhonov) path-of-solutions approach to hierarchical and bilevel variational inequalities~\cite{Kalashnikova}{, see also~\cite{Attouch1996},} and yields a general, if conservative, convergence rate. The second, presented in Section~\ref{sec:cent-eb}, exploits the error bound assumption (see Assumption~\ref{assum:penalty-eb} below) directly and yields a sharper rate.

% \begin{remark}
% Proposition~\ref{prop:reg-rate} below applies Lemma~\ref{lem:local-quad-growth} at the profile's own constant $S_0$, using only that the constraint violation along $\{\ab_\varepsilon^\star\}$ vanishes as $\varepsilon\to0$ -- valid for any admissible profile, with no need to enlarge $S_0$. For the quadratic profile of Example~\ref{ex:quad}, $S_0$ may in fact be taken arbitrarily large, so the lemma applies on \emph{every} bounded region; for the Huber profile of Example~\ref{ex:huber}, $S_0=M$ is fixed by the cap.
% \end{remark}

\subsection{Analysis based on regularized solutions}

To analyze the algorithm~\eqref{eq:algC}, we use the idea{ of}~\cite{Kalashnikova}, which {perturbs the lower-level operator $\Pb$ by a vanishing multiple of the upper-level operator $\Mb$} and reformulates the two-level problem as a single-level one.

Let us fix a positive parameter $\varepsilon$ and define $\Fb_{\varepsilon}(\ab):=\Pb(\ab)+\varepsilon\Mb(\ab)$. Consider the following regularized $\VI(\Ab,\Fb_{\varepsilon})$:
\begin{align}\label{eq:singleVI}
    \text{find } {\ab_\varepsilon^\star} \in \Ab \text{ such that }\cr
    \langle \Fb_{\varepsilon}({\ab_\varepsilon^\star}), \yb - {\ab_\varepsilon^\star} \rangle \geq 0,
    \quad \forall \yb \in \Ab.
    \end{align}
%\PL{In~\eqref{eq:singleVI} the solution was denoted $\ab^\star$, which is the v-GNE. I changed it to $\ab_\varepsilon^\star$, as in \cref{prop:reg-existence}.}
The next proposition connects solutions to $\VI(\Ab,\Fb_{\varepsilon})$ with those of the hierarchical VI~\eqref{eq:HVI}, and thus with the v-GNE of $\Gamma$.

\begin{proposition}\label{prop:reg-existence}
Let Assumptions~\ref{assum:convex} and~\ref{assum:strMon} hold. Then for {any} $\varepsilon>0$ the problem $\VI(\Ab,\Fb_{\varepsilon})$ has a unique solution $\ab_{\varepsilon}^\star$, and this solution converges to the unique v-GNE $\ab^\star$ of $\Gamma$ as $\varepsilon\to 0$. Moreover, for any two positive $\varepsilon_1$ and $\varepsilon_2$,
\[\|\ab^\star_{\varepsilon_1}-\ab^\star_{\varepsilon_2}\|\le\frac{|\varepsilon_2 - \varepsilon_1|}{\varepsilon_2\mu}\|\Mb(\ab_{\varepsilon_1}^\star)\|.\]
\end{proposition}

\begin{proof}
%\PL{I replaced the citation by this short proof, which needs no Slater condition. \eqref{eq:regM} and \eqref{eq:regbound} are reused in \cref{prop:reg-rate} and \cref{th:regularized}.}
{Since $\Pb$ is monotone and $\Mb$ is $\mu$-strongly monotone on $\Ab$, the continuous operator $\Fb_\varepsilon$ is $\varepsilon\mu$-strongly monotone on $\Ab$, so $\VI(\Ab,\Fb_\varepsilon)$ has a unique solution~\cite{FaccPang1}. Let $\zb\in S_L$. Since $\Psi$ is convex and $\Psi(\zb)=0$, we have $\Psi(\ab_\varepsilon^\star)\le\langle\Pb(\ab_\varepsilon^\star),\ab_\varepsilon^\star-\zb\rangle$. Combined with the VI inequality at $\yb=\zb$, this gives
\begin{align}\label{eq:regM}
0\le\Psi(\ab_\varepsilon^\star)\le\varepsilon\langle\Mb(\ab_\varepsilon^\star),\zb-\ab_\varepsilon^\star\rangle,\qquad\forall\zb\in S_L.
\end{align}
Taking $\zb=\ab^\star$ and using $\mu$-strong monotonicity of $\Mb$,
\begin{align}\label{eq:regbound}
\mu\|\ab_\varepsilon^\star-\ab^\star\|^2
&\le\langle\Mb(\ab_\varepsilon^\star)-\Mb(\ab^\star),\ab_\varepsilon^\star-\ab^\star\rangle\cr
&\le-\langle\Mb(\ab^\star),\ab_\varepsilon^\star-\ab^\star\rangle,
\end{align}
so that $\|\ab_\varepsilon^\star-\ab^\star\|\le\|\Mb(\ab^\star)\|/\mu$ for every $\varepsilon>0$. In particular, $\{\ab_\varepsilon^\star\}$ is bounded, so every sequence of positive numbers $\varepsilon_k\to0$ has a subsequence along which $\ab^\star_{\varepsilon_k}$ converges. As $\Mb$ is Lipschitz continuous, $\|\Mb(\ab_\varepsilon^\star)\|$ is bounded as well, and~\eqref{eq:regM} with $\zb=\ab^\star$ gives $\Psi(\ab_\varepsilon^\star)\to0$ as $\varepsilon\to0$. Now let $\bar\ab$ be a cluster point of $\{\ab_\varepsilon^\star\}$ as $\varepsilon\to0$, that is, $\ab^\star_{\varepsilon_k}\to\bar\ab$ for some sequence of positive numbers $\varepsilon_k\to0$. Since $\Ab$ is closed and $\Psi$ is continuous, $\bar\ab\in\Ab$ and $\Psi(\bar\ab)=0$, so $\bar\ab\in S_L$ by \cref{prop:pen-basic}. Dividing~\eqref{eq:regM} by $\varepsilon_k$ and passing to the limit gives $\langle\Mb(\bar\ab),\zb-\bar\ab\rangle\ge0$ for all $\zb\in S_L$, so $\bar\ab$ solves $\VI(S_L,\Mb)$. Since this VI has the unique solution $\ab^\star$, every cluster point equals $\ab^\star$. Hence $\ab_\varepsilon^\star\to\ab^\star$.}
Next, let $\ab_{\varepsilon_1}^\star$ and $\ab_{\varepsilon_2}^\star$ be the solutions to $\VI(\Ab,\Fb_{\varepsilon_1})$ and $\VI(\Ab,\Fb_{\varepsilon_2})$, respectively. Hence,
    \begin{align*}
        &\langle  \Pb(\ab_{\varepsilon_1}^\star)+\varepsilon_1\Mb(\ab_{\varepsilon_1}^\star), \ab_{\varepsilon_2}^\star - \ab_{\varepsilon_1}^\star \rangle \geq 0,\cr
        &\langle  \Pb(\ab_{\varepsilon_2}^\star)+\varepsilon_2\Mb(\ab_{\varepsilon_2}^\star),\ab_{\varepsilon_1}^\star - \ab_{\varepsilon_2}^\star \rangle \geq 0.
    \end{align*}
    Summing up the inequalities above and using monotonicity of $\Pb$ and $\mu$-strong monotonicity of $\Mb$ (Assumption~\ref{assum:strMon}), we conclude that
    \begin{align*}
        0
        &\le\langle\varepsilon_1\Mb(\ab_{\varepsilon_1}^\star)-\varepsilon_2\Mb(\ab_{\varepsilon_2}^\star), \ab_{\varepsilon_2}^\star - \ab_{\varepsilon_1}^\star\rangle \cr
        &= \langle\varepsilon_1\Mb(\ab_{\varepsilon_1}^\star)-\varepsilon_2\Mb(\ab_{\varepsilon_1}^\star), \ab_{\varepsilon_2}^\star - \ab_{\varepsilon_1}^\star\rangle \cr
        &\quad+ \langle\varepsilon_2\Mb(\ab_{\varepsilon_1}^\star)-\varepsilon_2\Mb(\ab_{\varepsilon_2}^\star), \ab_{\varepsilon_2}^\star - \ab_{\varepsilon_1}^\star\rangle\cr
        &\le |\varepsilon_1-\varepsilon_2|\|\Mb(\ab_{\varepsilon_1}^\star)\|\|\ab_{\varepsilon_2}^\star - \ab_{\varepsilon_1}^\star\| \cr
        &\quad- \varepsilon_2\mu\|\ab_{\varepsilon_2}^\star - \ab_{\varepsilon_1}^\star\|^2.
    \end{align*}
    Thus,
    \[\|\ab_{\varepsilon_2}^\star - \ab_{\varepsilon_1}^\star\| \le \frac{|\varepsilon_2 - \varepsilon_1|}{\varepsilon_2\mu}\|\Mb(\ab_{\varepsilon_1}^\star)\|.\]
\end{proof}
%\new{Note that \cref{prop:reg-existence} requires no error bound on the feasible set; \cref{assum:local-eb} is used below only to quantify the rate at which $\ab_\varepsilon^\star$ converges to $\ab^\star$.}

The following auxiliary result is an immediate consequence of admissibility (Definition~\ref{def:profile}).
\begin{lem}\label{lem:local-quad-growth}
Let $\psi$ be an admissible penalty profile with constants $S_0,c_\psi>0$ as in Definition~\ref{def:profile}. If $\ab\in{\R^D}$ satisfies $\|[\gb(\ab)]_+\|\le S_0$, then
$\Psi(\ab) \ge c_\psi\|[\gb(\ab)]_+\|^2$.
\end{lem}
\begin{proof}
Since $\|[\gb(\ab)]_+\|\le S_0$, every coordinate satisfies $[g_j(\ab)]_+\le\|[\gb(\ab)]_+\|\le S_0$, so Definition~\ref{def:profile} gives $\psi([g_j(\ab)]_+)\ge c_\psi[g_j(\ab)]_+\min([g_j(\ab)]_+,S_0)=c_\psi[g_j(\ab)]_+^2$ for every $j\in[m]$. Summing over $j$ and recalling~\eqref{eq:pen} gives the claim.
\end{proof}

The next proposition {complements \cref{prop:reg-existence} with a rate: it} characterizes the distance between $\ab_{\varepsilon}^\star$ (the solution to the regularized $\VI(\Ab,\Fb_{\varepsilon})$) and $\ab^\star$ (the v-GNE in $\Gamma$){ under the additional \cref{assum:local-eb}}.

\begin{proposition}\label{prop:reg-rate}
Let Assumptions~\ref{assum:convex},~\ref{assum:strMon}, and~\ref{assum:local-eb} hold. Then
    \[
\|\ab_\varepsilon^\star-\ab^\star\| = O(\sqrt{\varepsilon}) \quad \mbox{as }\varepsilon\to 0.
\]
\end{proposition}

\begin{proof}
Recall that $S_L = \Ab\cap \Cb$. For each $\varepsilon>0$, let $\sb_{\varepsilon}^\star = \Pi_{S_L}(\ab_{\varepsilon}^\star)$ and
$d_\varepsilon = \|\ab_{\varepsilon}^\star-\sb_{\varepsilon}^\star\|
= \operatorname{dist}(\ab_{\varepsilon}^\star,{S_L})$.
By Proposition~\ref{prop:reg-existence}, $\ab_{\varepsilon}^\star\to\ab^\star$ as $\varepsilon\to0$, so $\{\ab_{\varepsilon}^\star\}$ eventually lies in some fixed compact ball $B\subseteq\R^D$. {By \cref{assum:local-eb}, there exists} a constant $\kappa_B>0$ such that
\[
\operatorname{dist}(\ab,S_L) \le \kappa_B\|[\gb(\ab)]_+\|, \qquad \forall \ab\in B\cap\Ab.
\]
Moreover, since $\gb$ is continuous and $\gb(\ab^\star)\le\boldsymbol0$ (as $\ab^\star\in S_L$), $\|[\gb(\ab_\varepsilon^\star)]_+\|\to0$ as $\varepsilon\to0$. Hence, for all sufficiently small $\varepsilon$, $\ab_\varepsilon^\star\in B\cap\Ab$ and $\|[\gb(\ab_\varepsilon^\star)]_+\|\le S_0$, so that Lemma~\ref{lem:local-quad-growth} implies
\[
\Psi(\ab_\varepsilon^\star) = \sum_{j=1}^m\psi([g_j(\ab_\varepsilon^\star)]_+) \ge c_\psi\|[\gb(\ab_\varepsilon^\star)]_+\|^2.
\]
Combining the last inequalities above (applied at $\ab=\ab_\varepsilon^\star\in B\cap\Ab$), for all sufficiently small $\varepsilon$, we conclude that
\begin{align}\label{eq:regEB1}
 \Psi(\ab_{\varepsilon}^\star) \geq \frac{c_\psi}{\kappa_B^2} d_\varepsilon^2.
\end{align}

{By~\eqref{eq:regM} with $\zb=\sb_\varepsilon^\star$,}
\begin{align}\label{eq:regEB3}
\Psi(\ab_{\varepsilon}^\star)
\le
\varepsilon \langle \Mb(\ab_{\varepsilon}^\star),\sb_{\varepsilon}^\star-\ab_{\varepsilon}^\star\rangle
\le
\varepsilon \|\Mb(\ab_{\varepsilon}^\star)\|\, d_\varepsilon.
\end{align}
Since $\ab_{\varepsilon}^\star\to \ab^\star$ (Proposition~\ref{prop:reg-existence}) and $\Mb$ is continuous, there exists $\bar M>0$ such that
$\|\Mb(\ab_{\varepsilon}^\star)\| \le \bar M$ for sufficiently small $\varepsilon$. Therefore, by~\eqref{eq:regEB1} and~\eqref{eq:regEB3},
\[
\frac{c_\psi}{\kappa_B^2} d_\varepsilon^2
\le
\bar M \varepsilon d_\varepsilon.
\]
Thus, $d_\varepsilon \le \frac{\kappa_B^2}{c_\psi} \bar M\,\varepsilon$,
which proves $\operatorname{dist}(\ab_{\varepsilon}^\star,{S_L})=O(\varepsilon)$.

{By~\eqref{eq:regbound},
\begin{align}\label{eq:regdist1}
\mu\|\ab_{\varepsilon}^\star-\ab^\star\|^2
&\le-\langle \Mb(\ab^\star),\ab_{\varepsilon}^\star-\sb_{\varepsilon}^\star\rangle \cr
&\quad- \langle \Mb(\ab^\star),\sb_{\varepsilon}^\star-\ab^\star\rangle.
\end{align}}

Since $\sb_{\varepsilon}^\star\in {S_L}$ and $\ab^\star$ solves $\VI(S_L,\Mb)$, we have
$\langle \Mb(\ab^\star),\sb_{\varepsilon}^\star-\ab^\star\rangle \geq 0$.
Therefore, from~\eqref{eq:regdist1},
\[
\mu \|\ab_{\varepsilon}^\star-\ab^\star\|^2
\le
-\langle \Mb(\ab^\star),\ab_{\varepsilon}^\star-\sb_{\varepsilon}^\star\rangle
\le
\|\Mb(\ab^\star)\|\,\|\ab_{\varepsilon}^\star-\sb_{\varepsilon}^\star\|.
\]
Hence, $\mu \|\ab_{\varepsilon}^\star-\ab^\star\|^2
\le
\|\Mb(\ab^\star)\|\, \operatorname{dist}(\ab_{\varepsilon}^\star,{S_L})$.
Since $\operatorname{dist}(\ab_{\varepsilon}^\star,{S_L})=O(\varepsilon)$, we conclude that $\|\ab_{\varepsilon}^\star-\ab^\star\| = O(\sqrt{\varepsilon})$.
\end{proof}
%\PL{This proof uses the convergence $\ab_\varepsilon^\star\to\ab^\star$ of \cref{prop:reg-existence} only for boundedness of $\{\ab_\varepsilon^\star\}$, which~\eqref{eq:regbound} gives for every $\varepsilon>0$, and for $\|[\gb(\ab_\varepsilon^\star)]_+\|\le S_0$ for small $\varepsilon$, which also follows from $\Psi(\ab_\varepsilon^\star)\to0$. Convergence also follows from the $O(\sqrt\varepsilon)$ rate proved here. So, technically, the convergence part of \cref{prop:reg-existence} can be dropped, and the sensitivity bound could be derived in the proof of \cref{th:regularized}, unless we want to keep \cref{prop:reg-existence} as a result of interest in its own right. What do you think?}

Next we prove the following convergence result for the centralized algorithm~\eqref{eq:algC}.

\begin{theorem}\label{th:regularized}
Let Assumptions~\ref{assum:convex}--\ref{assum:Lip} hold. Let the stepsize $\gamma$ satisfy $0 < \gamma < \frac{2}{L_{\Pb}}$, let $\eta>0$, and let the regularization parameter be set up as $\varepsilon_t = \frac{E}{(t+1)^e}$ for some $E>0${ and \(e \in (0,1)\)}. Then the sequence generated by~\eqref{eq:algC} converges to the unique v-GNE $\ab^\star$ in $\Gamma$ and
\[\|\ab(t) - \ab^{\star}\|^2= O\left(\frac{1}{t^{\min\{e, 2-2e\}}}\right).\]
%\PL{We need to restrict the range for \(e\) otherwise this rate does not make sense.}
\end{theorem}
\begin{proof}
Set $\tau_t = \frac{\eta}{\gamma}\varepsilon_t$, and denote by $\ab^{\star}_{t} = \ab^{\star}_{\tau_t}$ the unique solution to $\VI(\Ab,\Pb+\tau_t \Mb)$. Since $\tau_t = \frac{\eta E/\gamma}{(t+1)^e}$, Propositions~\ref{prop:reg-existence} and~\ref{prop:reg-rate} apply with $\varepsilon$ replaced by $\tau_t$.

By the projection identity for VI solutions, written in the form of~\eqref{eq:algC}, we have
\begin{align*}
\ab^{\star}_{t}
&=
\Pi_{\Ab}\bigl[
\ab^{\star}_{t} - \gamma(\Pb(\ab^{\star}_{t})+\tau_t \Mb(\ab^{\star}_{t}))
\bigr]\cr
&=
\Pi_{\Ab}\bigl[
\ab^{\star}_{t} - \gamma\Pb(\ab^{\star}_{t})-\eta\varepsilon_t \Mb(\ab^{\star}_{t})
\bigr].
\end{align*}
Define $\Fb_t(\ab) = \Fb_{\tau_t}(\ab) = \Pb(\ab) + \tau_t \Mb(\ab)$, so that $\gamma \Fb_t(\ab)=\gamma\Pb(\ab)+\eta\varepsilon_t\Mb(\ab)$ matches~\eqref{eq:algC}.%\PL{Here $F_t$ and $F_{\tau_t}$ are not bold, while~\eqref{eq:singleVI} uses $\Fb_\varepsilon$. Also $F$ denotes the lifted operator in \cref{rem:lVI}. Should we write $\Fb_t$ here?}
 Then by the nonexpansiveness of $\Pi_{\Ab}$, we obtain
\begin{align}\label{eq:regTh1}
\|\ab(t+1) - \ab^{\star}_{t}\|^2
&\le
\left\|
(\ab(t) - \ab^{\star}_{t})
- \gamma\bigl( \Fb_t(\ab(t)) - \Fb_t(\ab^{\star}_{t}) \bigr)
\right\|^2\cr
& = \|\ab(t) - \ab^{\star}_{t}\|^2 \cr
&\quad- 2\gamma\langle \Fb_t(\ab(t)) - \Fb_t(\ab^{\star}_{t}), \ab(t) - \ab^{\star}_{t} \rangle \cr
&\quad+ \gamma^2\| \Fb_t(\ab(t)) - \Fb_t(\ab^{\star}_{t})\|^2.
\end{align}
Using cocoercivity of $\Pb$ over $\Ab$ (which holds since {$\Psi$ is convex, }$\Pb=\nabla\Psi$ and $\Pb$ is $L_{\Pb}$-Lipschitz { on the open convex set $U\supseteq\Ab$ by Assumption~\ref{assum:Lip}}, see \cite[Thm.~3.1]{PerezArosVilches2021}), %\PL{Lipschitz continuity on $\Ab$ is not enough to invoke Baillon-Haddad, unless we assume $\Ab$ has nonempty interior (which is too much to ask). I changed \cref{assum:Lip} to require Lipschitz continuity on an open convex neighborhood $U\supseteq\Ab$. We can then apply Baillon-Haddad on $U$ and restrict the inequality to $\Ab$. I also adjusted the sufficient conditions after the assumption. The stepsize restrictions use the Lipschitz constant on $U$ now. \cite[Thm.~3.3]{BauschkeCombettes2010} requires $\Psi\in C^2$, which fails for our profiles at $g_j=0$, so I cite \cite[Thm.~3.1]{PerezArosVilches2021}, which only needs differentiability.}
\[
\langle \Pb(\ab(t))-\Pb(\ab^\star_t), \ab(t) - \ab^{\star}_{t}\rangle
\ge \frac{1}{L_{\Pb}}\|\Pb(\ab(t))-\Pb(\ab^\star_t)\|^2,
\]
together with Lipschitz continuity and $\mu$-strong monotonicity of $\Mb$ (Assumption~\ref{assum:strMon}), and the weighted Young inequality $\|\ub+\vb\|^2\le(1+\delta)\|\ub\|^2+(1+\delta^{-1})\|\vb\|^2$ with $1+\delta=\frac{2}{\gamma L_{\Pb}}$ (so that $\delta>0$ by $\gamma<\frac{2}{L_{\Pb}}$ and $1+\delta^{-1}=K_\gamma:=\frac{2}{2-\gamma L_{\Pb}}$), we obtain from~\eqref{eq:regTh1}
\begin{align}\label{eq:regTh1b}
   \|\ab(t+1) - \ab^{\star}_{t}\|^2
\le
&\|\ab(t) - \ab^{\star}_{t}\|^2 -2\gamma\tau_t\mu\|\ab(t) - \ab^{\star}_{t}\|^2 \cr
&- \frac{2\gamma}{L_{\Pb}}\|\Pb(\ab(t))-\Pb(\ab^\star_t)\|^2 \cr
&+ K_\gamma\gamma^2\tau_t^2{L_{\Mb}^2}\|\ab(t) - \ab^{\star}_{t}\|^2 \cr
&+ \frac{2\gamma}{L_{\Pb}}\|\Pb(\ab(t))-\Pb(\ab^\star_t)\|^2 \cr
=&(1-2\gamma\tau_t\mu+K_\gamma\gamma^2\tau_t^2L_{\Mb}^2)\|\ab(t) - \ab^{\star}_{t}\|^2\cr
=&(1-2\eta\varepsilon_t\mu+K_\gamma\eta^2\varepsilon_t^2L_{\Mb}^2)\|\ab(t) - \ab^{\star}_{t}\|^2,
\end{align}
where we used $(1+\delta)\gamma^2=\frac{2\gamma}{L_{\Pb}}$, and in the last equality we substituted $\gamma\tau_t=\eta\varepsilon_t$.

Next,
\[
\|\ab(t+1) - \ab^{\star}_{t+1}\|
\le
\|\ab(t+1) - \ab^{\star}_{t}\|
+
\|\ab^{\star}_{t+1} - \ab^{\star}_{t}\|.
\]

By sensitivity of strongly monotone VIs (see Proposition~\ref{prop:reg-existence}, applied to the sequence $\tau_t${ with $\varepsilon_1=\tau_{t+1}$ and $\varepsilon_2=\tau_t$}), {with $C_0=\sup_{\tau>0}\|\Mb(\ab^\star_\tau)\|/\mu$, which is finite by~\eqref{eq:regbound} and Lipschitz continuity of $\Mb$,}
\[
\|\ab^{\star}_{t+1} - \ab^{\star}_{t}\|
\le
C_0 \frac{|\tau_{t+1} - \tau_t|}{\tau_t}
=
C_0 \frac{|\varepsilon_{t+1} - \varepsilon_t|}{\varepsilon_t}\old{,}\new{.}
\]
Thus, taking into account the last three inequalities, we conclude that
\begin{align*}
\|\ab(t+1) - \ab^{\star}_{t+1}\|
\le{}&
 {\Bigl(1-\tfrac{\eta\mu}{2}\varepsilon_t\Bigr)}\|\ab(t) - \ab^{\star}_{t}\|\cr
&+
C_0 \frac{|\varepsilon_{t+1} - \varepsilon_t|}{\varepsilon_t},
\end{align*}
{for all $t$ large enough that $\varepsilon_t\le\min\{\mu/(K_\gamma\eta L_{\Mb}^2),1/(\eta\mu)\}$, since then the factor in~\eqref{eq:regTh1b} is at most $1-\eta\mu\varepsilon_t$ and $\sqrt{1-\eta\mu\varepsilon_t}\le1-\frac{\eta\mu}{2}\varepsilon_t$. Moreover, $|\varepsilon_{t+1}-\varepsilon_t|/\varepsilon_t=1-(1-\frac1{t+2})^e\le\frac1{t+2}$.}
%\PL{\(c\) depended on \(t\). I fixed the rate to \(\eta\mu/2\) and removed the square root, so that \cref{lem:chung}(i) applies.}
{Since $\varepsilon_t=E/(t+1)^e$ with $e\in(0,1)$}, {\cref{lem:chung}(i) with $p=e$, $q=1$, and $a=\eta\mu E/2$ implies}
$\|\ab(t) - \ab^{\star}_{t}\| = O\left(\frac{1}{t^{1-e}}\right)$. Taking into account Proposition~\ref{prop:reg-rate} (applied to $\tau_t$, of the same order as $\varepsilon_t$){ and $\|\ab(t)-\ab^\star\|\le\|\ab(t)-\ab^\star_t\|+\|\ab^\star_t-\ab^\star\|$}, we conclude that $\|\ab(t) - \ab^{\star}\| = O\left(\frac{1}{t^{\min\{e/2, 1-e\}}}\right)${, and squaring gives the claim}.%\PL{I replaced the citation of Chung's lemma by \cref{lem:chung}(i) in the appendix. It needs a lower bound on \(\varepsilon_t\), so \(O(1/t^e)\) was the wrong direction.}
\end{proof}
\begin{corollary}\label{cor:regularized-optimal}
Optimizing the rate in Theorem~\ref{th:regularized} above, we obtain that under the choice $e=2/3$,
\[\|\ab(t) - \ab^{\star}\|^2 = O\left(\frac{1}{t^{2/3}}\right).\]
\end{corollary}

%\PL{--------------------------- Checked until here -----------------------}

\subsection{Error Bound Condition}\label{sec:cent-eb}

We now show that strengthening the local error bound of Assumption~\ref{assum:local-eb} to the global one posed below, i.e., requiring a single constant $\kappa$ valid on all of $\Ab$, yields a sharper rate than Theorem~\ref{th:regularized} and Corollary~\ref{cor:regularized-optimal} above: in particular, a rate of $O(1/t)$ becomes available.
To this end, we pose the following additional assumption on the constraint set.
\begin{assumption}\label{assum:penalty-eb}
The constraint set in~\eqref{eq:HVI} admits an error bound on $\Ab$, i.e., there exists
$\kappa>0$ such that
\[
\operatorname{dist}(\ab,S_L)
\le
\kappa \|[\gb(\ab)]_+\|,
\qquad \forall \ab\in \Ab.
\]
\end{assumption}
{Assumption~\ref{assum:penalty-eb} is the global counterpart of \cref{assum:local-eb}: it implies \cref{assum:local-eb} with $\kappa_B=\kappa$ for every bounded $B$, but the converse fails in general, even under Slater's condition. For instance, let $\Ab=\R^2$ and $g(\ab)=\|\ab\|-a_1-1$, so that $\Ab\cap \Cb=\{\ab:\,a_2^2\le 2a_1+1\}$ and $\ab=\boldsymbol 0$ is a Slater point; hence \cref{assum:local-eb} holds by \cref{lem:slater-eb}. Along $\ab=[x,\sqrt{2x+1}+1]$, however, $\operatorname{dist}(\ab,\Ab\cap \Cb)\to1$ while $[g(\ab)]_+\to0$ as $x\to\infty$, so Assumption~\ref{assum:penalty-eb} fails.}

\begin{remark}\label{rem:AffineConst}
Assumption~\ref{assum:penalty-eb} is satisfied, in particular, when $\gb$ is affine and the individual action set $\Ab$ is polyhedral,
namely $\gb(\ab)=K\ab-\boldsymbol k$ and $\Ab = \{\ab:\, H\ab- \boldsymbol h\le \boldsymbol{0}\}$.
Indeed, %in this case $\Psi(\ab)=\frac12\|[B\ab-\bb]_+\|^2$.
by Hoffman's error bound (see, for example,~\cite{Hoffman}), there exists $\kappa>0$ such that
\begin{align*}
\operatorname{dist}(\ab,S_L) &= \operatorname{dist}(\ab,\Ab\cap \Cb)
\le
\kappa\bigl(\|[H\ab-\boldsymbol h]_+\|+\|[K\ab-\boldsymbol k]_+\|\bigr) \cr&= \kappa\|[K\ab-\boldsymbol k]_+\|,
\end{align*}
where we used $\ab\in \Ab$.

{Another sufficient condition for Assumption~\ref{assum:penalty-eb} is \cref{assum:local-eb} together with boundedness of $\Ab$: taking $B=\Ab$ gives Assumption~\ref{assum:penalty-eb} with $\kappa=\kappa_{\Ab}$. In particular, this is the case under Slater's condition with bounded $\Ab$, where \cref{lem:slater-eb} gives $\kappa=\sup_{\ab\in\Ab}\|\ab-\bar\ab\|/\bar s$, with $\bar s$ as in \cref{lem:slater-eb}.}
\end{remark}

% Admissibility of $\psi$ (Definition~\ref{def:profile}) already quantifies $\Psi$ in terms of the constraint violation $\|[\gb(\ab)]_+\|$, so no further assumption on $\psi$ itself is needed to convert the error bound of Assumption~\ref{assum:penalty-eb} into a bound on $\Psi$ over all of $\Ab$. (The stronger, purely quadratic growth condition below remains useful in its own right -- it is what the distributed convergence analysis of Section~\ref{sec:distributed} relies on.)
% \begin{assumption}\label{assum:quad-growth}
% There exists $c_\psi>0$ such that $\psi(s)\ge c_\psi s^2$ for all $s\ge0$.
% \end{assumption}
% \begin{remark}
% The quadratic profile of Example~\ref{ex:quad} satisfies Assumption~\ref{assum:quad-growth} on all of $[0,\infty)$ with $c_\psi=\frac12$; the Huber profile of Example~\ref{ex:huber} does not, since it grows only linearly beyond its cap $M$. Assumption~\ref{assum:quad-growth} is \emph{not} needed for Lemma~\ref{lem:centralized-eb} or Theorem~\ref{th:centralized} below, which rely only on Definition~\ref{def:profile}; it is retained here for the distributed convergence analysis.
% \end{remark}
% \TT{Comment on the Assumption}

Converting the per-coordinate growth of Definition~\ref{def:profile} into a bound on $\Psi$ in terms of the \emph{norm} $\|[\gb(\ab)]_+\|$ relies on the following elementary fact.
\begin{lem}\label{lem:norm-hybrid}
Let $s_1,\ldots,s_m\ge0$, $S_0>0$, and $r=\bigl(\sum_{j=1}^m s_j^2\bigr)^{1/2}$. Then
$\sum_{j=1}^m s_j\min(s_j,S_0) \;\ge\; r\min(r,S_0)$.
\end{lem}
\begin{proof}
Let $\mathcal J_1=\{j:s_j\le S_0\}$, $\mathcal J_2=\{j:s_j>S_0\}$, and set $T_1=\sum_{\mathcal J_1}s_j^2$, $T_2=\sum_{\mathcal J_2}s_j^2$ (so $r^2=T_1+T_2$) and $T_3=\sum_{\mathcal J_2}s_j$; the left-hand side equals $T_1+S_0T_3$, and $T_3\ge\sqrt{T_2}$ since $(\sum_{\mathcal J_2}s_j)^2\ge\sum_{\mathcal J_2}s_j^2$ for nonnegative terms.

If $r\le S_0$, every $s_j\le r\le S_0$, so $\mathcal J_2=\emptyset$ and the left-hand side equals $T_1=r^2=r\min(r,S_0)$.

If $r>S_0$, it suffices, using $T_3\ge\sqrt{T_2}$, to show $T_1+S_0\sqrt{T_2}\ge S_0\sqrt{A+T_2}$; since both sides are nonnegative, this is equivalent to $T_1^2+2T_1S_0\sqrt{T_2}\ge S_0^2T_1$, i.e., (for $T_1>0$) to $T_1+2S_0\sqrt{T_2}\ge S_0^2$, with the case $T_1=0$ verified directly below. As a function of $T_1\ge\max(0,S_0^2-T_2)$ (the constraint $T_1+T_2> S_0^2$), the left-hand side is increasing in $T_1$, so it suffices to check it at $T_1=\max(0,S_0^2-T_2)$: if $T_2\le S_0^2$, this gives $(S_0^2-T_2)+2S_0\sqrt{T_2}=S_0^2+\sqrt{T_2}(2S_0-\sqrt{T_2})\ge S_0^2$ (as $\sqrt{T_2}\le S_0<2S_0$); if $T_2>S_0^2$, this gives $2S_0\sqrt{T_2}>2S_0\cdot S_0=2S_0^2>S_0^2$. Either way, $T_1+2S_0\sqrt{T_2}\ge S_0^2$, completing the proof.
\end{proof}

The following lemma, which relies only on Definition~\ref{def:profile}{, the convexity in Assumption~\ref{assum:convex}(ii),} and the error bound of Assumption~\ref{assum:penalty-eb}, is the key tool used to establish this sharper convergence rate for~\eqref{eq:algC}.

\begin{lem}\label{lem:centralized-eb}
Let $\zb\in S_L$, let $\psi$ be an admissible penalty profile with constants $S_0,c_\psi>0$ as in Definition~\ref{def:profile}, and set $D_0:=\kappa S_0$. Then, under Assumptions~\ref{assum:convex}(ii) and~\ref{assum:penalty-eb}, for any $\ab\in\Ab$ and given $d=\operatorname{dist}(\ab,S_L)$,
\[
\langle \Pb(\ab),\ab-\zb\rangle
\ge
\mu_{\Pb}\, d\min(d,D_0), \qquad \mu_{\Pb}=\frac{c_\psi}{\kappa^2}.
\]
\end{lem}
\begin{proof}
Let $r=\|[\gb(\ab)]_+\|$ and $d=\operatorname{dist}(\ab,S_L)$. By Definition~\ref{def:profile} applied coordinatewise and Lemma~\ref{lem:norm-hybrid}, recalling~\eqref{eq:pen},
\begin{align*}
\Psi(\ab) = \sum_{j=1}^m\psi([g_j(\ab)]_+) &\ge c_\psi\sum_{j=1}^m[g_j(\ab)]_+\min([g_j(\ab)]_+,S_0) \cr
&\ge c_\psi\, r\min(r,S_0).
\end{align*}
By Assumption~\ref{assum:penalty-eb}, $r\ge d/\kappa$, and since $s\mapsto s\min(s,S_0)$ is nondecreasing on $[0,\infty)$,
\begin{align*}
\Psi(\ab) &\ge c_\psi\cdot\frac d\kappa\min\Bigl(\frac d\kappa,S_0\Bigr) \cr
&= \frac{c_\psi}{\kappa^2}\,d\min(d,\kappa S_0) = \mu_{\Pb}\, d\min(d,D_0).
\end{align*}
Since $\Psi$ is convex and $\Psi(\zb)=0$ for $\zb\in S_L$,
\[
0=\Psi(\zb)\ge \Psi(\ab)+\langle \Pb(\ab),\zb-\ab\rangle,
\]
i.e., $\langle \Pb(\ab),\ab-\zb\rangle\ge \Psi(\ab)\ge \mu_{\Pb}\, d\min(d,D_0)$.
\end{proof}

\begin{theorem}\label{th:centralized}
Let Assumptions~\ref{assum:convex},~\ref{assum:strMon},~\ref{assum:Lip}, and~\ref{assum:penalty-eb} hold. Let the regularization parameter be $\varepsilon_t=\frac{E}{(t+1)^e}$ for some $E>0$, $e\in(0,1]$, let the stepsize satisfy $\gamma\in\left(0,\frac{2}{L_{\Pb}}\right)$, and let $\eta>0$. Then the sequence $\{\ab(t)\}$ generated by~\eqref{eq:algC} converges to the unique v-GNE $\ab^\star$ in $\Gamma$. Moreover,
\[
\|\ab(t)-\ab^\star\|^2
=
\begin{cases}
&O\left(\frac1{t^e}\right), \mbox{ if } e\in(0,1),\\
&O\left(\frac1{t}\right), \mbox{ if } e=1 \mbox{ and } E>\frac{1}{\eta\mu}.
\end{cases}
\]
\end{theorem}

\begin{proof}
Set $\db_t=\ab(t)-\ab^\star$. Since $\ab^\star\in S_L=\operatorname{SOL}(\Ab,\Pb)$, we have
\[
\ab^\star =\Pi_{\Ab}[\ab^\star-\gamma \Pb(\ab^\star)].
\]
Using nonexpansiveness of the projection,
\[
\|\ab(t+1)-\ab^\star\|^2
\le
\left\|\db_t-\gamma(\Pb(\ab(t))-\Pb(\ab^\star))
-\eta\varepsilon_t \Mb(\ab(t))\right\|^2 .
\]
Let $\pb_t=\Pb(\ab(t))-\Pb(\ab^\star)$. Fix $\delta=\frac{2-\gamma L_{\Pb}}{2\gamma L_{\Pb}}$, which is positive since $\gamma<\frac{2}{L_{\Pb}}$, and set $\bar K_\gamma:=1+\delta^{-1}=\frac{2+\gamma L_{\Pb}}{2-\gamma L_{\Pb}}$. Then, since $2\gamma\eta\varepsilon_t \langle\pb_t,\Mb(\ab(t))\rangle \le \delta\gamma^2\|\pb_t\|^2 + \delta^{-1}\eta^2\varepsilon_t^2\|\Mb(\ab(t))\|^2$, we obtain
\begin{align}\label{eq:cent1}
    \|\ab(t+1)-\ab^\star\|^2
&\le
\|\db_t\|^2
-2\gamma\langle \pb_t,\db_t\rangle
-2\eta\varepsilon_t\langle \Mb(\ab(t)),\db_t\rangle
\cr
&\quad
+(1+\delta)\gamma^2\|\pb_t\|^2
+\bar K_\gamma\eta^2\varepsilon_t^2\|\Mb(\ab(t))\|^2.
\end{align}
Since $\Pb=\nabla\Psi$ is the gradient of a convex function with an $L_{\Pb}$-Lipschitz gradient{ on the open convex set $U\supseteq\Ab$} (Assumption~\ref{assum:Lip}), it is cocoercive{ on $\Ab$ by \cite[Thm.~3.1]{PerezArosVilches2021}}, i.e., $\|\pb_t\|^2\le L_{\Pb}\langle \pb_t,\db_t\rangle$. Thus, since $(1+\delta)L_{\Pb}\gamma^2=\frac{\gamma(2+\gamma L_{\Pb})}{2}$, setting $\bar\gamma:=\frac{\gamma(2-\gamma L_{\Pb})}{2}>0$ gives
\begin{align*}
-2\gamma\langle \pb_t,\db_t\rangle
+(1+\delta)\gamma^2\|\pb_t\|^2
&\le -2\gamma\langle \pb_t,\db_t\rangle
+(1+\delta)L_{\Pb}\gamma^2\langle \pb_t,\db_t\rangle\cr
&=
-\bar\gamma\langle \pb_t,\db_t\rangle.
\end{align*}
Thus,~\eqref{eq:cent1} implies
\begin{align}\label{eq:cent2}
    \|\ab(t+1)-\ab^\star\|^2
&\le
\|\db_t\|^2
-\bar\gamma\langle \pb_t,\db_t\rangle
-2\eta\varepsilon_t\langle \Mb(\ab(t)),\db_t\rangle\cr
&\quad+\bar K_\gamma\eta^2\varepsilon_t^2\|\Mb(\ab(t))\|^2 .
\end{align}
Let $d_t=\operatorname{dist}(\ab(t),S_L)$. As $\ab^\star\in S_L$, $\Pb(\ab^\star) = 0$. Then Lemma~\ref{lem:centralized-eb} implies
$\langle \pb_t,\db_t\rangle= \langle \Pb(\ab(t)),\ab(t)-\ab^\star\rangle\ge
\mu_{\Pb}\, d_t\min(d_t,D_0)$,
so that
\begin{align}\label{eq:cent3a}
-\bar\gamma\langle \pb_t,\db_t\rangle
\le
-\bar\gamma\mu_{\Pb}\, d_t\min(d_t,D_0).
\end{align}
Next decompose
\[
\langle \Mb(\ab(t)),\db_t\rangle
=
\langle \Mb(\ab(t))-\Mb(\ab^\star),\db_t\rangle
+
\langle \Mb(\ab^\star),\db_t\rangle.
\]
By $\mu$-strong monotonicity of $\Mb$ (Assumption~\ref{assum:strMon}),
\[
\langle \Mb(\ab(t))-\Mb(\ab^\star),\db_t\rangle
\ge
\mu\|\db_t\|^2.
\]
Since $\ab^\star$ solves $\VI(S_L,\Mb)$,
$\langle \Mb(\ab^\star),\Pi_{S_L}(\ab(t))-\ab^\star\rangle\ge 0$.
Therefore,
\[
\begin{aligned}
\langle \Mb(\ab^\star),\db_t\rangle
&=
\langle \Mb(\ab^\star),\ab(t)-\Pi_{S_L}(\ab(t))\rangle
\\
&\quad+
\langle \Mb(\ab^\star),\Pi_{S_L}(\ab(t))-\ab^\star\rangle
\\
&\ge
-\|\Mb(\ab^\star)\|\, d_t.
\end{aligned}
\]
Thus,
\begin{align}\label{eq:cent3b}
\langle \Mb(\ab(t)),\db_t\rangle
\ge
\mu\|\db_t\|^2
-
\|\Mb(\ab^\star)\|\, d_t.
\end{align}
Combining~\eqref{eq:cent3a} and~\eqref{eq:cent3b} in~\eqref{eq:cent2} yields
\begin{align}\label{eq:cent3}
\|\ab(t+1)-\ab^\star\|^2
&\le
\|\db_t\|^2
-\bar\gamma\mu_{\Pb}\, d_t\min(d_t,D_0)
\cr
&\quad
-2\eta\mu\varepsilon_t\|\db_t\|^2
\cr
&\quad+
2\eta\varepsilon_t\|\Mb(\ab^\star)\|\, d_t
\cr
&\quad
+\bar K_\gamma\eta^2\varepsilon_t^2\|\Mb(\ab(t))\|^2 .
\end{align}
We bound the two $d_t$-dependent terms in~\eqref{eq:cent3} jointly. Let $\bar\varepsilon=\frac{\bar\gamma\mu_{\Pb}D_0}{4\eta\|\Mb(\ab^\star)\|}$ and $C_3=\frac{2\eta^2}{\bar\gamma\mu_{\Pb}}\|\Mb(\ab^\star)\|^2$. Then for all $t$ with $\varepsilon_t\le\bar\varepsilon$,
\begin{align}\label{eq:cent-hybrid}
-\bar\gamma\mu_{\Pb}\, d_t\min(d_t,D_0) + 2\eta\varepsilon_t\|\Mb(\ab^\star)\|\, d_t \le C_3\varepsilon_t^2.
\end{align}
Indeed, if $d_t\le D_0$, then $\min(d_t,D_0)=d_t$ and~\eqref{eq:cent-hybrid} is the standard Young's-inequality bound $2\eta\varepsilon_t\|\Mb(\ab^\star)\|d_t\le C_3\varepsilon_t^2+\bar\gamma\mu_{\Pb}d_t^2$, obtained by completing the square in $d_t$. If $d_t>D_0$, then $\min(d_t,D_0)=D_0$ and, since $\varepsilon_t\le\bar\varepsilon$ gives $2\eta\varepsilon_t\|\Mb(\ab^\star)\|\le\frac{\bar\gamma\mu_{\Pb}D_0}{2}$, the left-hand side of~\eqref{eq:cent-hybrid} equals
\begin{align*}
d_t\bigl(2\eta\varepsilon_t\|\Mb(\ab^\star)\|-\bar\gamma\mu_{\Pb}D_0\bigr)
&\le
d_t\Bigl(\tfrac{\bar\gamma\mu_{\Pb}D_0}{2}-\bar\gamma\mu_{\Pb}D_0\Bigr)\cr
&=
-\tfrac{\bar\gamma\mu_{\Pb}D_0}{2}\,d_t
\;\le\;
0
\;\le\;
C_3\varepsilon_t^2.
\end{align*}
Since $\varepsilon_t\to0$, $\varepsilon_t\le\bar\varepsilon$ for all sufficiently large $t$. Combining~\eqref{eq:cent3} with~\eqref{eq:cent-hybrid} then gives, for such $t$,
\begin{align}\label{eq:cent4}
\|\ab(t+1)-\ab^\star\|^2
\le{}&
\|\db_t\|^2
-2\eta\mu\varepsilon_t\|\db_t\|^2\cr
&+\bar K_\gamma\eta^2\varepsilon_t^2\|\Mb(\ab(t))\|^2
+C_3\varepsilon_t^2 .
\end{align}
Since $\Mb$ is Lipschitz continuous (Assumption~\ref{assum:strMon}),
\begin{align*}
\bar K_\gamma\|\Mb(\ab(t))\|^2
\le C_1\|\db_t\|^2+C_2,
\end{align*}
where $C_1=2\bar K_\gamma{L_{\Mb}^2}$ and $C_2=2\bar K_\gamma\|\Mb(\ab^\star)\|^2$. Hence, inserting the last inequality into~\eqref{eq:cent4} gives, for all sufficiently large $t$,
\begin{align*}
\|\ab(t+1)-\ab^\star\|^2
\le{}&
\left(1-2\eta\mu\varepsilon_t+\eta^2C_1\varepsilon_t^2\right)
\|\ab(t)-\ab^\star\|^2\cr
&+
(C_3+\eta^2C_2)\varepsilon_t^2.
\end{align*}
Due to the choice of $\varepsilon_t$, for all sufficiently large $t$,
$
2\eta\mu\varepsilon_t-\eta^2C_1\varepsilon_t^2
\ge
\eta\mu\varepsilon_t$,
so that, for such $t$,
\[
\|\ab(t+1)-\ab^\star\|^2
\le
(1-\eta\mu\varepsilon_t)\|\ab(t)-\ab^\star\|^2
+
(C_3+\eta^2C_2)\varepsilon_t^2.
\]
{By \cref{lem:chung} with $p=e$, $q=2e$, and $a=\eta\mu E$, we obtain convergence, by item (iii) if $e=1$ and $\eta\mu E\le1$, and}
\[
\|\ab(t)-\ab^\star\|^2
=
\begin{cases}
&O\left(\frac1{t^e}\right), \mbox{ if } e\in(0,1),\\
&O\left(\frac1{t}\right), \mbox{ if } e=1 \mbox{ and } E>\frac{1}{\eta\mu}.
\end{cases}
\]
\end{proof}

\section{Distributed Hierarchical Regularization Scheme}\label{sec:distributed}

In this section, we develop a distributed algorithmic framework for solving the hierarchical variational inequality~\eqref{eq:HVI}. Since the mappings \(\Mb\) and \(\Pb\) depend on the global action profile, the centralized scheme of Section~\ref{sec:centralized} would generally require global information, namely, access of every agent to the full joint action $\ab(t)$. To overcome this limitation, we adopt a consensus-based distributed approach in which each agent maintains local estimates of the actions of all other agents and communicates only with its neighbors over a communication network.

Specifically, each agent $i$
stores a local estimate of the global action profile
$\ab = [\ab^1,\ldots,\ab^N]$ which is denoted by
\[
\xb_{i}
=
[\xb_{i}^1,\ldots,\xb_{i}^N]\in\R^D.
\]
Here, $\xb_{i}^j$ represents agent $i$'s estimate of the action of agent $j$. In particular, $\xb_{i}^i=\ab^i$.
Stacking all local estimates yields
\[
\xb
=
[\xb_{1},\ldots,\xb_{N}
]
\in\R^{ND}.
\]

The agents communicate over an undirected graph
$\mathcal G=([N],\mathcal E)$,
where \((i,j)\in\mathcal E\) if agents \(i\) and \(j\) can exchange information. For each agent \(i\in[N]\), we denote by $\mathcal N_i=\{j\in[N]\mid (i,j)\in\mathcal E\}
$
the set of its neighbors. We assume that \(i\in\mathcal N_i\), that is, each agent is considered a neighbor of itself. Throughout the paper, we impose the following standard assumption~\cite{olshevsky2009convergence}.
\begin{assumption}\label{assum:comm}
The communication graph \(\mathcal G\) is connected. Moreover, the associated mixing matrix
$W=[w_{ij}]_{i,j=1}^N$ defines the weights on the undirected arcs such that $w_{ij}>0$  
if and only if the agents $i$ and $j$ can exchange information (are neighbors). $W$
is symmetric and doubly stochastic, i.e.,
$W=W^\top$, $W\boldsymbol 1_N=\boldsymbol 1_N$, and $\boldsymbol 1_N^\top W=\boldsymbol 1_N^\top$.
\end{assumption}
The assumption above implies, in particular, that the second largest singular value of $W$ satisfies
\begin{align}\label{eq:sigma}
 \sigma=\|W-\bar J\|\in[0,1), \qquad \bar J=\frac{1}{N}\boldsymbol 1_N\boldsymbol 1_N^\top.    
\end{align}
Assumption~\ref{assum:comm} guarantees that repeated applications of the mixing matrix asymptotically drive the agents' local estimates to consensus. The parameter \(\sigma\) quantifies the connectivity of the communication graph and determines the rate of consensus formation~\cite{olshevsky2009convergence}.
%Specifically, we consider the time-varying operator 
% \begin{align}\label{eq:Ft}
% \Fb_t=\Pb+\varepsilon_t\Mb,    
% \end{align}
% where
% $\varepsilon_t>0$ is chosen such that $\varepsilon_t\to0$ as time runs.
% Due to this setting, the influence of the upper-level operator $\Mb$ gradually vanishes while preserving the hierarchical selection property.
For the current analysis, we assume the following condition regarding the feasible sets and smoothness of the operators.
\begin{assumption}\label{assum:box}
Each set $\Ab^i$, $i\in[N]$, is bounded. Fix a compact box
\[
\Qb^i\supseteq \Ab^i,\quad i\in[N],\qquad
\Qb:=\Qb^1\times\cdots\times \Qb^N\supseteq\Ab,
\]
known to all agents (e.g., $\Qb^i=[-\rho,\rho]^{d_i}$ for $\rho>0$ large enough that $\Ab^i\subseteq \Qb^i$). Moreover, $\Pb$ and $\Mb$ are Lipschitz continuous over $\Qb$ with constants $L_{\Pb}^{\Qb},L_{\Mb}^{\Qb}>0$, respectively:
\begin{align*}
\|\Pb(\ab)-\Pb(\bb)\|&\le L_{\Pb}^{\Qb}\|\ab-\bb\|,
\cr
\|\Mb(\ab)-\Mb(\bb)\|&\le L_{\Mb}^{\Qb}\|\ab-\bb\|,
\qquad\forall \ab,\bb\in \Qb.
\end{align*}
\end{assumption}

Using the notation above, we define the steps of the distributed procedure in Algorithm~\ref{alg:1}. The algorithm consists of three main steps. First, each agent $i$ exchanges its current estimates with its neighbors and computes a local weighted average $\zb_i$ using the corresponding row of the mixing matrix $W$. Second, each agent $i$ performs a projected gradient step with respect to the local coordinate of the regularized operator, where the lower-level part $\Pb$ is calculated at the local weighted average $\zb_i$, whereas the upper-level part $\Mb^i$ is calculated at the joint action estimate $\xb_i$ and influences only the agent's own action update. Finally, the resulting iterate of the local action update is projected onto the local feasible set $\Ab^i$ for the owned coordinate, and onto the shared public box $\Qb^j$ (see Assumption~\ref{assum:box}) for every other, non-owned coordinate $j\neq i$.
Importantly, the proposed algorithm operates entirely in the primal space and therefore does not require introducing, communicating, or updating dual variables; the boxes $\Qb^j$ are common knowledge and disclose no private information about the sets $\Ab^j$.

\begin{algorithm}[t]
\caption{Distributed hierarchical regularization scheme}
\label{alg:1}
\begin{algorithmic}[1]
\State \textbf{Initialization:} Each agent chooses $\xb_{i}(0)\in\boldsymbol{\Omega}^i=\Qb^1\times\cdots\times\Qb^{i-1}\times\Ab^i\times\Qb^{i+1}\times\cdots\times\Qb^N$. The stepsizes $\gamma,\eta>0$ and the regularization sequence $\{\varepsilon_t\}\subseteq(0,\infty)$ are common to all agents.
\For{$t=0,1,2,\ldots$}
\State Each agent $i\in[N]$ obtains the estimates $\xb_{j}(t)$ from all neighbors $j\in\mathcal N_i$ and computes \[ \zb_{i}(t) = \sum_{j=1}^N w_{ij}\xb_{j}(t). \]
\State Each agent $i\in[N]$ updates its local estimate as
\begin{align*}
\xb^i_i(t+1) &= \Pi_{\Ab^i} \Bigl[ \zb^i_{i}(t) - \gamma\Pb^i(\zb_{i}(t)) - \eta\varepsilon_t\Mb^i(\xb_{i}(t)) \Bigr], \cr
 \xb_{i}^j(t+1)&=\Pi_{\Qb^j}\Bigl[\zb_{i}^j(t)-\gamma \Pb^j(\zb_{i}(t))\Bigr], \qquad j\neq i.
\end{align*}
\EndFor
\end{algorithmic}
\end{algorithm}
\begin{remark}
Assumption~\ref{assum:box} is of exactly the same nature as Assumption~\ref{assum:Lip}, only stated on the compact set $\Qb\supseteq\Ab$ instead of {the open neighborhood $U$ of $\Ab$}, and is satisfied under the same sufficient conditions listed after Assumption~\ref{assum:Lip} with {$U$} replaced by {an open convex neighborhood of }$\Qb$ throughout. {Since $\Qb$ is a box with nonempty interior, cocoercivity of $\Pb$ on $\operatorname{int}\Qb$ follows from \cite[Thm.~3.1]{PerezArosVilches2021} and extends to $\Qb$ by continuity of $\Pb$.}  Assumption~\ref{assum:box} imposes no control whatsoever on $\Pb$ or $\Mb$ outside the fixed, bounded, and explicitly constructed set $\Qb$. Agent $i$ need not know the private set $\Ab^j$ of any other agent $j$. Assumption~\ref{assum:box} allows for keeping $\zb_i(t)$ inside a fixed, compact, and explicitly known set on which $\Pb$ is well controlled.
\end{remark}

For each agent $i$, let us define the selection matrix $R_i\in\R^{d_i\times D}$: 
\begin{align}\label{eq:R}
    \setlength{\arraycolsep}{2.5pt}
R_i
=
\begin{bmatrix}
0_{d_i\times d_1}
&
\cdots
&
0_{d_i\times d_{i-1}}
&
I_{d_i}
&
0_{d_i\times d_{i+1}}
&
\cdots
&
0_{d_i\times d_N}
\end{bmatrix},
\end{align}
where, as before, \(D=\sum_{i=1}^N d_i\).

% Since $\Pb=\nabla\Psi$ with $\Psi$ convex on all of $\R^D$, $\Psi\ge0$, and $\Psi(a^\star)=0=\phi_{\min}$ (so that $\Pb(a^\star)=\boldsymbol 0$, $a^\star$ being a global, unconstrained minimizer of $\Psi$), the standard Baillon--Haddad argument, applied along the convex set $Q$ (large enough, per Assumption~\ref{assum:box}, to also contain the argument's auxiliary points), turns $L_{\Pb}^{\Qb}$-Lipschitz continuity of $\Pb$ on $Q$ into the co-coercivity inequality
% \[
% \left\langle \Pb(\ab)-\Pb(a^\star),\ab-a^\star\right\rangle
% \ge
% \frac{1}{L_{\Pb}^{\Qb}}\|\Pb(\ab)-\Pb(a^\star)\|^2,
% \qquad\forall \ab\in \Qb,
% \]
% which is the form of Assumption~\ref{assum:box} used in the proof of Theorem~\ref{th:one-level-constant-step} below.

We further define the sets of local joint action estimates:
% \begin{align}\label{eq:R}
% R  = \operatorname{diag}(R_1,\ldots,R_N)
% \in\R^{D\times ND},
% \end{align}
\[
\boldsymbol{\Omega}^i
=
\Qb^1
\times\cdots\times
\Qb^{i-1}
\times\Ab^i
\times
\Qb^{i+1}
\times\cdots\times
\Qb^N
\]
and $\boldsymbol{\Omega}=
\prod_{i=1}^{N}\boldsymbol{\Omega}^i$.
% Since $\Ab^i\subseteq \Qb^i$ for every $i\in[N]$, we have $\Ab\subseteq \Qb$, and every $\xb\in\boldsymbol\Omega$ satisfies $x_i^j\in \Qb^j$ for all $i,j\in[N]$ (with $j=i$ using $\Ab^i\subseteq \Qb^i$). Consequently, for $\boldsymbol z=\boldsymbol W\boldsymbol x$, each component
% $
% z_i^j=\sum_{l=1}^N w_{il}x_l^j
% $
% is a convex combination of points of the convex set $\Qb^j$, hence $z_i^j\in \Qb^j$. Thus
% \[
% z_i\in \Qb,\qquad\forall i\in[N],\ \forall\xb\in\boldsymbol\Omega.
% \]
% Together with $a(\xb)\in\Ab\subseteq \Qb$, this places both arguments of $\Pb$ needed below inside $\Qb$, so that Assumption~\ref{assum:box} -- rather than Assumption~\ref{assum:Lip} -- applies in the proof of Lemma~\ref{lem:one-level-residual-control}.
% For
% \[
% \xb
% =
% \operatorname{col}(x_1,\ldots,x_N)
% \in\R^{ND},
% \]
% define the owned-action vector
% \[
% a(\xb)
% :=
% \operatorname{col}
% \bigl(
% R_1x_1,\ldots,R_Nx_N
% \bigr).
% \]
%Notice that
%\[
%\xb\in\boldsymbol{\Omega}
%\quad\Longrightarrow\quad
%a(\xb)\in\mathcal{A}.
%\]
Next, for $\xb = [\xb_1,\ldots,\xb_N]\in\R^{ND}$, we define, with a slight abuse of notation,
\begin{align}\label{eq:one-level-operatorP}
\Pb(\xb)
=
[
\Pb(\xb_1),\ldots,\Pb(\xb_N)
]
\end{align}
and
\[
\Mb_R(\xb)
=
[
R_1^\top \Mb^1(\xb_1),\ldots,
R_N^\top \Mb^N(\xb_N)
].
\]
Moreover, define
$\Wb
=
W\otimes I_D$. Then the compact form of the distributed hierarchical regularization scheme is given by
\begin{align}
\label{eq:one-level-algorithm}
\zb(t) &= \Wb\xb(t),\cr
\xb(t+1)&=
\Pi_{\boldsymbol{\Omega}}
\left[
\zb(t)
-
\gamma\Pb(\zb(t))
-\eta\varepsilon_t
\Mb_R(\xb(t))
\right].
\end{align}

\subsection{Convergence Analysis}\label{sec:conv_distr}

We first formulate two auxiliary results (Lemmas~\ref{lem:one-level-owned-disagreement} and~\ref{lem:one-level-residual-control}) that are used to conduct the convergence analysis of Algorithm~\ref{alg:1}. Their proofs are deferred to the Appendix.

% -------------------------------------------------------------------------
% Preliminary lemmas
% -------------------------------------------------------------------------
Next, for
$\xb=[\xb_1,\ldots,\xb_N]
\in\R^{ND}$, we define the owned-action vector
\[
\ab(\xb)
=
[
R_1\xb_1,\ldots,R_N\xb_N
], 
\]
where $R_i$, $i\in[N]$, is defined in~\eqref{eq:R}.
Notice that, if $\xb\in\boldsymbol{\Omega}$, then $\ab(\xb)\in\Ab$.
Moreover, we define
\[
\Pi_\perp
=
(I_N-\bar J)\otimes I_D, \qquad \xb^\perp
=
\Pi_\perp\xb,\]
as well as the consensual lower-level solution set:
\[
\boldsymbol{S}_L
=
\left\{
\boldsymbol 1_N\otimes \ab:\ab\in S_L
\right\}.
\]

\begin{lem}%[Owned-action disagreement]
\label{lem:one-level-owned-disagreement}
For every $\xb\in\R^{ND}$,
$
\left\|\xb-\boldsymbol 1_N\otimes \ab(\xb)
\right\|^2
\le
(N+1)\|\xb^\perp\|^2$.
\end{lem}

\begin{lem}%[Residual control]
\label{lem:one-level-residual-control}
Let $\psi$ be an admissible penalty profile as in Definition~\ref{def:profile}. Under Assumptions~\ref{assum:convex}{(ii)},~\ref{assum:local-eb},~\ref{assum:comm}, and~\ref{assum:box}, define, for $\xb\in\boldsymbol{\Omega}$,
$
\zb
=
\Wb\xb$
and
$
\mathcal{R}(\xb)
=
\|\xb^\perp\|^2
+
\|\Pb(\zb)\|^2$.

Then there exists \(C_R>0\), independent of \(\xb\), such that
\[
\operatorname{dist}
\left(
\xb,
\boldsymbol{S}_L
\right)^2
\leq
C_R\mathcal{R}(\xb),
\qquad
\forall\,\xb\in\boldsymbol{\Omega}.
\]
Moreover, there exist constants \(K_R,K_a>0\), independent of \(\xb\) and \(\sb\), such that
\[
\begin{aligned}
\|\xb-\boldsymbol 1_N\otimes \sb\|^2
&\leq
K_R\mathcal{R}(\xb)
+
K_a\|\ab(\xb)-\sb\|^2,
\\
&\qquad
\forall\,\xb\in\boldsymbol{\Omega},\ \forall\,\sb\in S_L.
\end{aligned}
\]
\end{lem}

% -------------------------------------------------------------------------
% Main theorem
% -------------------------------------------------------------------------
With the results of Lemmas~\ref{lem:one-level-owned-disagreement} and~\ref{lem:one-level-residual-control} in hand, we are now ready to present the main convergence result for Algorithm~\ref{alg:1}.
Throughout the remainder of this section, the regularization parameter in Algorithm~\ref{alg:1} is taken to be of the form
\[
\varepsilon_t=\frac{E}{(t+1)^e},\qquad E>0,\qquad e\in(0,1].
\]
%so that $\varepsilon_t\to0$, $\sum\varepsilon_t=\infty$, and $\sum\varepsilon_t^2<\infty$.

\begin{theorem}
\label{th:one-level-constant-step}
Let Assumptions~\ref{assum:convex},~\ref{assum:strMon},~\ref{assum:local-eb},~\ref{assum:comm}, and~\ref{assum:box} hold. Let \(\ab^\star\) be the unique solution of $\VI(S_L,\Mb)$,
and define
$\xb^\star
=
\boldsymbol 1_N\otimes \ab^\star$.
Then, for any $\gamma\in(0,\frac{2}{L_{\Pb}^{\Qb}})
$ and any $\eta>0$,
the sequence generated by~\eqref{eq:one-level-algorithm} satisfies
$
\xb(t)\to\xb^\star
$ as $t\to\infty$.        Moreover,
\[
\|\xb(t)-\xb^\star\|^2
=
\begin{cases}
&O\left(\frac1{t^e}\right), \mbox{ if } e\in(0,1),\\
&O\left(\frac1{t}\right), \mbox{ if } e=1 \mbox{ and } E>\frac{1}{c}.
\end{cases}
\]  
where \(c>0\) is a constant, given explicitly in the proof, that depends only on the problem and the algorithm parameters.
\end{theorem}
\begin{proof}
For the sake of brevity, set $\xb_t=\xb(t)$ with blocks $\xb_{i,t}=\xb_i(t)$, $\zb_t=\zb(t)$, $\db_t=
\xb_t-\xb^\star$, and
$\ab_t=\ab(\xb_t)$.
Define
\begin{align}\label{eq:calR}
\mathcal{R}_t=
\|\xb_t^\perp\|^2
+
\|\Pb(\zb_t)\|^2.
\end{align}
Since \(\ab^\star\in S_L\),
$\Pb(\ab^\star)=\boldsymbol{0}$.
Therefore,
$\Pb(\xb^\star)=\boldsymbol{0}$.
Define
\[
\ub_t
=
\zb_t-\xb^\star
-
\gamma\Pb(\zb_t).
\]
Since $\xb_t\in\boldsymbol\Omega$, we have $\zb_{i,t}\in \Qb$ for all $i\in[N]$ (due to double stochasticity of $W$)
and $\xb^\star=\boldsymbol 1_N\otimes \ab^\star\in \Qb\times\cdots\times\Qb$ (as $\ab^\star\in S_L\subseteq\Ab\subseteq \Qb$). Thus the cocoercivity inequality implied by Assumption~\ref{assum:box}{ and \cite[Thm.~3.1]{PerezArosVilches2021} on $\operatorname{int}\Qb$, extended to $\Qb$ by continuity of $\Pb$,} applies to the pair $(\zb_t,\xb^\star)$, resulting in $\langle \Pb(\zb_t),\zb_t-\xb^\star\rangle = \langle \Pb(\zb_t)-\Pb(\xb^\star),\zb_t-\xb^\star\rangle\ge \frac{1}{L_{\Pb}^{\Qb}}\|\Pb(\zb_t)-\Pb(\xb^\star)\|^2=\frac{1}{L_{\Pb}^{\Qb}}\|\Pb(\zb_t)\|^2$. Hence,
\begin{align}\label{eq:1}
\|\ub_t\|^2
&=
\|\zb_t-\xb^\star\|^2
-
2\gamma
\left\langle
\Pb(\zb_t),
\zb_t-\xb^\star
\right\rangle
+
\gamma^2
\|\Pb(\zb_t)\|^2
\cr
&\leq
\|\zb_t-\xb^\star\|^2
-
\gamma
\left(
\frac{2}{L_{\Pb}^{\Qb}}-\gamma
\right)
\|\Pb(\zb_t)\|^2.
\end{align}
Moreover, since $\xb^\star=\boldsymbol 1_N\otimes \ab^\star$ is already consensual, $\Pi_\perp\xb^\star=\boldsymbol 0$, and therefore $\Pi_\perp\db_t=\Pi_\perp\xb_t=\xb_t^\perp$. Recall $\bar J=\frac1N\boldsymbol 1_N\boldsymbol 1_N^\top$. By double stochasticity of $W$ (Assumption~\ref{assum:comm}), one has $W\bar J=\bar JW=\bar J$, so that
$\Wb(\bar J\otimes I_D)=\bar J\otimes I_D$, $\Wb\Pi_\perp=(W-\bar J)\otimes I_D$, and $\Wb\xb^\star=\xb^\star$.
Consequently,
\[
\zb_t-\xb^\star
=
\Wb\db_t
=
(\bar J\otimes I_D)\db_t
+
\bigl((W-\bar J)\otimes I_D\bigr)\xb_t^\perp.
\]
The two terms on the right are orthogonal, because $\bar J(W-\bar J)=\bar JW-\bar J^2=\bar J-\bar J=0$. Moreover, since $\bar J\otimes I_D$ and $\Pi_\perp$ are complementary orthogonal projections, $\|\db_t\|^2=\|(\bar J\otimes I_D)\db_t\|^2+\|\xb_t^\perp\|^2$. Combining these relations with the spectral bound $\|W-\bar J\|=\sigma\in[0,1)$ furnished by Assumption~\ref{assum:comm} (see~\eqref{eq:sigma}) gives
\begin{align*}
\|\zb_t-\xb^\star\|^2
&=
\|(\bar J\otimes I_D)\db_t\|^2
+
\bigl\|\bigl((W-\bar J)\otimes I_D\bigr)\xb_t^\perp\bigr\|^2\cr
&\leq
\|\db_t\|^2
-
(1-\sigma^2)\|\xb_t^\perp\|^2.
\end{align*}
Hence,~\eqref{eq:1} and the preceding inequalities yield (see also~\eqref{eq:calR})
\begin{align}\label{eq:u}
\|\ub_t\|^2
\leq
\|\db_t\|^2
-
c_L\mathcal{R}_t,
\end{align}
where
$
c_L
=
\min
\left\{
1-\sigma^2,\,
\gamma
\left(
\frac{2}{L_{\Pb}^{\Qb}}-\gamma
\right)
\right\}
>
0$,
which is strictly positive precisely because $0<\gamma<2/L_{\Pb}^{\Qb}$.

Since \(\xb^\star\in\boldsymbol{\Omega}\),
nonexpansiveness of the projection in~\eqref{eq:one-level-algorithm}
gives
\begin{align}\label{eq:fin1}
\|\xb_{t+1}-\xb^\star\|^2
&\leq
\left\|
\ub_t
-
\eta\varepsilon_t
\Mb_R(\xb_t)
\right\|^2
\cr
&=
\|\ub_t\|^2
-
2\eta\varepsilon_t
\left\langle
\Mb_R(\xb_t),
\ub_t
\right\rangle
\cr
&\quad+
\eta^2\varepsilon_t^2
\|\Mb_R(\xb_t)\|^2.
\end{align}
We now bound the middle and last terms of this expansion. The constants, obtained from the compactness furnished by Assumption~\ref{assum:box}, will be used throughout. Since $\Ab$ is compact, the diameter $\Delta=\sup_{\ab,\bb\in\Ab}\|\ab-\bb\|$ 
is finite. Moreover, the maximum
$
M_{\max}=\max_{\ab\in \Ab}\|\Mb(\ab)\|
$
is attained and finite. Because $\xb_t\in\boldsymbol\Omega$ for every $t$, we have $\ab_t,\ab^\star\in\Ab$ and $\xb_{i,t}\in \Qb$ for every $i$, so that, throughout the proof,
$
\|\ab_t-\ab^\star\|\leq\Delta$ and $\|\Mb^i(\xb_{i,t})\|\leq\|\Mb(\xb_{i,t})\|\leq M_{\max}$ for all $i\in[N]$ and all $t$. 
%\smallskip

\noindent\textbf{Estimate of the middle term $\left\langle
\Mb_R(\xb_t),
\ub_t
\right\rangle$.}
Since
\[
\ub_t
=
\db_t
+
(\zb_t-\xb_t)
-
\gamma\Pb(\zb_t),
\]
we have
\[
\begin{aligned}
\left\langle
\Mb_R(\xb_t),
\ub_t
\right\rangle
&=
\underbrace{\left\langle
\Mb_R(\xb_t),
\db_t
\right\rangle}_{\text{(a)}}
\\
&\quad+
\underbrace{\left\langle
\Mb_R(\xb_t),
\zb_t-\xb_t
\right\rangle
-
\gamma
\left\langle
\Mb_R(\xb_t),
\Pb(\zb_t)
\right\rangle}_{\text{(b)}}.
\end{aligned}
\]

\smallskip
\noindent\emph{Term (a).}
By definition of \(\Mb_R\) and of the selection matrices $R_i$,
\begin{align}\label{eq:M_R}
\left\langle
\Mb_R(\xb_t),
\db_t
\right\rangle
=
\langle \Mb(\ab_t),\ab_t-\ab^\star\rangle
+
\zeta_t,
\end{align}
with $
\zeta_t
=
\sum_{i=1}^{N}
\left\langle
\Mb^i(\xb_{i,t})-\Mb^i(\ab_t),
\ab_t^i-\ab^{\star,i}
\right\rangle
$,
where $\ab_t^i=R_i\xb_{i,t}$ and $\ab^{\star,i}=R_i\ab^\star$ denote the $i$-th blocks of $\ab_t=\ab(\xb_t)$ and $\ab^\star$. By the Cauchy--Schwarz inequality, using $\sum_{i=1}^N\|\ab_t^i-\ab^{\star,i}\|^2=\|\ab_t-\ab^\star\|^2\leq\Delta^2$,
\[
|\zeta_t|
\leq
\Delta
\left(
\sum_{i=1}^{N}
\|\Mb^i(\xb_{i,t})-\Mb^i(\ab_t)\|^2
\right)^{1/2}.
\]
Since $\xb_{i,t},\ab_t\in \Qb$, Assumption~\ref{assum:box}  gives $\|\Mb^i(\xb_{i,t})-\Mb^i(\ab_t)\|\le L_{\Mb}^{\Qb}\|\xb_{i,t}-\ab_t\|$ for every $i$. Thus,
\[
\begin{aligned}
\left(
\sum_{i=1}^{N}
\|\Mb^i(\xb_{i,t})-\Mb^i(\ab_t)\|^2
\right)^{1/2}
&\leq
L_{\Mb}^{\Qb}
\left\|
\xb_t-\boldsymbol 1_N\otimes \ab_t
\right\|
\\
&\leq
L_{\Mb}^{\Qb}\sqrt{N+1}\,\|\xb_t^\perp\|,
\end{aligned}
\]
where in the last inequality we used Lemma~\ref{lem:one-level-owned-disagreement}.
Combining the last two relations and using $\|\xb_t^\perp\|\leq\sqrt{\mathcal{R}_t}$ (see again~\eqref{eq:calR}), we conclude that
\begin{align}\label{eq:E}
|\zeta_t|
\leq
\Delta L_{\Mb}^{\Qb}\sqrt{N+1}\,\sqrt{\mathcal{R}_t}.
\end{align}
For the remaining term, strong monotonicity of $\Mb$ (Assumption~\ref{assum:strMon}) together with the variational-inequality optimality of $\ab^\star$ for $\VI(S_L,\Mb)$, applied at $\sb_t=\Pi_{S_L}(\ab_t)$, i.e., $\langle \Mb(\ab^\star),\sb_t-\ab^\star\rangle\geq0$, gives
\begin{align}\label{eq:Mb}
\langle \Mb(\ab_t),\ab_t-\ab^\star\rangle
&=
\langle \Mb(\ab_t)-\Mb(\ab^\star),\ab_t-\ab^\star\rangle
\cr
&\quad+
\langle \Mb(\ab^\star),\ab_t-\ab^\star\rangle
\cr
&\geq
\mu\|\ab_t-\ab^\star\|^2
-
M_{\max}\bar d_t,
\end{align}
where $\bar d_t=\operatorname{dist}(\xb_t,\boldsymbol{S}_L)$, we used $\|\Mb(\ab^\star)\|\le M_{\max}$, and the last inequality uses $\|\ab_t-\sb_t\|=\operatorname{dist}(\ab_t,S_L)\leq \bar d_t$ (since $\boldsymbol 1_N\otimes \sb_t\in\boldsymbol S_L$ and $\|\ab_t-\sb_t\|\le\|\xb_t-\boldsymbol 1_N\otimes \sb_t\|$). By Lemma~\ref{lem:one-level-residual-control}, $\bar d_t\leq\sqrt{C_R}\sqrt{\mathcal{R}_t}$. Combining this last inequality with \eqref{eq:E} and \eqref{eq:Mb} and inserting it into \eqref{eq:M_R}, we obtain
\[\left\langle
\Mb_R(\xb_t),
\db_t
\right\rangle
\geq
\mu\|\ab_t-\ab^\star\|^2
-
C_a\sqrt{\mathcal{R}_t},
\]
where $C_a=M_{\max}\sqrt{C_R}+\Delta L_{\Mb}^{\Qb}\sqrt{N+1}$.

\smallskip
\noindent\emph{Term (b).}
Since $\zb_t-\xb_t=(\Wb-I_{ND})\xb_t = \bigl((W-I_N)\otimes I_D\bigr)\xb_t^\perp$, we have
\[
\|\zb_t-\xb_t\|
\leq
C_W\sqrt{\mathcal{R}_t},
\qquad
C_W=\|W-I_N\|.
\]
%and, directly from the definition of $\mathcal{R}_t$, $\|\Pb(\zb_t)\|\leq\sqrt{\mathcal{R}_t}$. 
Moreover, since
$\Mb_R(\xb_t)=[R_1^\top \Mb^1(\xb_{1,t}),\ldots,R_N^\top \Mb^N(\xb_{N,t})]$,
we have $
\|\Mb_R(\xb_t)\|^2
=
\sum_{i=1}^N\|\Mb^i(\xb_{i,t})\|^2
\leq
NM_{\max}^2$ and, thus, the uniform bound $\|\Mb_R(\xb_t)\|\leq\bar M_R=\sqrt N M_{\max}$. By the Cauchy--Schwarz inequality,
\[
\begin{aligned}
|\text{(b)}|
&\leq
\|\Mb_R(\xb_t)\|
\bigl(
\|\zb_t-\xb_t\|
+
\gamma\|\Pb(\zb_t)\|
\bigr)
\leq
C_b\sqrt{\mathcal{R}_t},
\end{aligned}\]
where $C_b=\bar M_R(C_W+\gamma)$.

\smallskip
\noindent\textbf{Combined estimate.}
Adding the bounds for (a) and (b), with $C_\ast=C_a+C_b$,
\begin{align}\label{eq:middle}
\left\langle
\Mb_R(\xb_t),
\ub_t
\right\rangle
\geq
\mu\|\ab_t-\ab^\star\|^2
-
C_\ast\sqrt{\mathcal{R}_t}.
\end{align}
Substituting \eqref{eq:u}, \eqref{eq:middle}, and $\|\Mb_R(\xb_t)\|^2\leq\bar M_R^2$ into the  estimate~\eqref{eq:fin1} gives
\[
\begin{aligned}
\|\xb_{t+1}-\xb^\star\|^2
\leq{}&
\|\db_t\|^2
-
c_L\mathcal{R}_t
-
2\eta\mu\varepsilon_t\|\ab_t-\ab^\star\|^2
\\
&+
2\eta C_\ast\varepsilon_t\sqrt{\mathcal{R}_t}
+
\eta^2\bar M_R^2\varepsilon_t^2.
\end{aligned}
\]
By Young's inequality, $2\eta C_\ast\varepsilon_t\sqrt{\mathcal{R}_t}\leq\frac{c_L}{2}\mathcal{R}_t+\frac{2\eta^2C_\ast^2}{c_L}\varepsilon_t^2$, so that, for every $t\geq0$,
\begin{align}\label{eq:fin2}
\|\xb_{t+1}-\xb^\star\|^2
\leq{}&
\|\db_t\|^2
-
\frac{c_L}{2}\mathcal{R}_t
-
2\eta\mu\varepsilon_t\|\ab_t-\ab^\star\|^2+
C_4\varepsilon_t^2,
\end{align}
where $C_4=\frac{2\eta^2C_\ast^2}{c_L}+\eta^2\bar M_R^2$.

By Lemma~\ref{lem:one-level-residual-control}, $\|\db_t\|^2\leq K_R\mathcal{R}_t+K_a\|\ab_t-\ab^\star\|^2$ for all $t\geq0$. Since $\varepsilon_t=E/(t+1)^e\leq E$ for every sufficiently large $t\geq0$, setting
$
c
=
\min
\left\{
\frac{c_L}{2K_RE},\ \frac{2\eta\mu}{K_a}
\right\}
$
ensures $c\varepsilon_tK_R\leq c_L/2$ and $cK_a\leq2\eta\mu$ for every $t\geq0$, and therefore
\[
\begin{aligned}
c\varepsilon_t\|\db_t\|^2
&\leq
c\varepsilon_tK_R\mathcal{R}_t
+
c\varepsilon_tK_a\|\ab_t-\ab^\star\|^2
\\
&\leq
\frac{c_L}{2}\mathcal{R}_t
+
2\eta\mu\varepsilon_t\|\ab_t-\ab^\star\|^2,
\end{aligned}
\]
implying $-\frac{c_L}{2}\mathcal{R}_t\le -c\varepsilon_t\|\db_t\|^2+2\eta\mu\varepsilon_t\|\ab_t-\ab^\star\|^2$. 
Consequently,~\eqref{eq:fin2} implies, for every sufficiently large $t\geq0$,
\[
\|\xb_{t+1}-\xb^\star\|^2
\leq
(1-c\varepsilon_t)
\|\xb_t-\xb^\star\|^2
+
C_4\varepsilon_t^2.
\]
{\Cref{lem:chung} with $p=e$, $q=2e$, and $a=cE$, applied to the scalar recursion above, gives}%\PL{$\sum_t\varepsilon_t^2<\infty$ fails for $e\le1/2$, and it is not needed. For $e=1$ and $cE\le1$, convergence follows from \cref{lem:chung}(iii).} 
$\|\xb_t-\xb^\star\|^2\to0$, and moreover $\|\xb_t-\xb^\star\|^2=O(\varepsilon_t)=O(1/t^e)$ whenever $e\in(0,1)$, while for $e=1$ one obtains the faster rate $\|\xb_t-\xb^\star\|^2=O(1/t)$ provided $cE>1$. This proves the convergence and rate claims of the theorem.
\end{proof}

\section{Simulation Results}\label{sec:simulations}

We evaluate the proposed hierarchical VI reformulation and its centralized and distributed algorithms (Sections~\ref{sec:centralized} and~\ref{sec:distributed}) on randomly generated strongly monotone GNEP instances constructed so that a ground-truth v-GNE $(\ab^\star,\bla^\star)$ satisfying the KKT conditions~\eqref{eq:KKT} is known by design, which allows us to directly measure the distance $\|\ab(t)-\ab^\star\|$ to the true equilibrium. We consider three settings: centralized GNE seeking with nonlinear shared constraints, and distributed GNE seeking with nonlinear and with linear shared constraints, respectively. In each setting, we report results averaged over $20$ random instances of problems of the same size, with strong monotonicity constant $\mu=0.1$. All tests were run on a MacBook Pro with an Apple M5 Max chip.

\subsection{Centralized, nonlinear shared constraints}
\label{sec:sim-centralized-nonlinear}
% example_hvi_beats_all.py
The test instances are obtained by drawing a target equilibrium $\ab^\star\in\R^D$ with i.i.d.\ standard Gaussian entries
and then constructing the GNEP around this target equilibrium with convex quadratic costs, local box constraints, and nonlinear shared constraints.
Following the construction suggested in~\cite{BT26c}, each player $i\in[N]$ has the cost 
\begin{align}\label{eq:sim-player}
J^i(\ab)=\tfrac12\ab^\top H^i\ab+(c^i)^\top\ab^i
\end{align}
where $c^i\in\R^{d_i}$ and the symmetric matrix $H^i$ has $(i,j)$ block $G_{ij}$ and $(j,i)$ block $G_{ij}^\top$ for all $j\in[N]$, and zero blocks elsewhere, with $G_{ij}\in\R^{d_i\times d_j}$ the blocks of the matrix
\[
    G=V^\top V+S-S^\top+\left(\mu-\lambda_{\min}(V^\top V)\right)I_D.
\]
Here $V\in\R^{D\times D}$ is upper triangular and $S\in\R^{D\times D}$ is strictly lower block triangular with respect to the partition of $\ab$ into $\ab^1,\dots,\ab^N$, both with i.i.d.\ standard Gaussian nonzero entries. By construction, the pseudo-gradient is the affine mapping $\Mb(\ab)=G\ab+c$, with $c=[c^1,\dots,c^N]$, and since $\frac12(G+G^\top)=V^\top V+(\mu-\lambda_{\min}(V^\top V))I$ has minimum eigenvalue equal to $\mu$, $\Mb$ is $\mu$-strongly monotone and $\|G\|$-Lipschitz continuous, as required by Assumption~\ref{assum:strMon}. The players are coupled through $m=n_e+n_q$ shared constraints $\gb(\ab)\leq\boldsymbol 0$, namely $n_e$ log-sum-exp and $n_q$ convex quadratic constraints
\begin{align}\label{eq:sim-coupling}
g_s(\ab)&=\log\sum_{k=1}^{2}\exp\left(\alpha_{s,k}^\top\ab+\beta_{s,k}\right)-r_s\leq 0 && s\in[n_e]\cr
g_{n_e+s}(\ab)&=\ab^\top P_s\ab+p_s^\top\ab-r_{n_e+s}\leq 0 && s\in[n_q]
\end{align}
where $\alpha_{s,k},p_s\in\R^D$ and $\beta_{s,k}\in\R$ have i.i.d.\ Gaussian entries (with variance $1/D$ for $\alpha_{s,k}$ and $p_s$), and $P_s=B_s^\top B_s\succeq0$, with $B_s\in\R^{D\times D}$ having i.i.d.\ Gaussian entries with variance $1/D$. A subset $\mathcal I_a\subseteq[m]$ of $m_a$ constraints is selected at random to be active at $\ab^\star$: the offsets $r_j$ are chosen such that $g_j(\ab^\star)=0$ for $j\in\mathcal I_a$ and $g_j(\ab^\star)=-\delta_j$ for $j\notin\mathcal I_a$, with slack $\delta_j$ uniformly distributed in $[0.5,1.5]$. The local action sets are boxes $\Ab^i=\{\ab^i:\ \underline\ab^i\leq\ab^i\leq\overline\ab^i\}$, whose bounds are obtained by subtracting from and adding to each entry of $\ab^{\star,i}$ numbers uniformly distributed in $[0.5,1.5]$, so that $\ab^\star$ lies in the interior of $\Ab$ and $\EuScript N_{\Ab}(\ab^\star)=\{\boldsymbol 0\}$. The multipliers are set to $\lambda^\star_j$ uniformly distributed in $[0.5,1.5]$ for $j\in\mathcal I_a$ and $\lambda^\star_j=0$ otherwise. Finally, the linear cost terms are chosen as
\[
    c=-G\ab^\star-\nabla\gb(\ab^\star)^\top\bla^\star
\]
so that the KKT conditions~\eqref{eq:KKT} hold at $(\ab^\star,\bla^\star)$ with strict complementarity. Since the game is strongly monotone with continuous costs over the defined compact, convex, nonempty constraint set, $\ab^\star$ is the unique v-GNE of the generated game. Such a random problem construction is implemented in \texttt{NashOpt}~\cite{Bem26}.

\begin{figure}[t]
\centering
\includegraphics[width=0.9\linewidth]{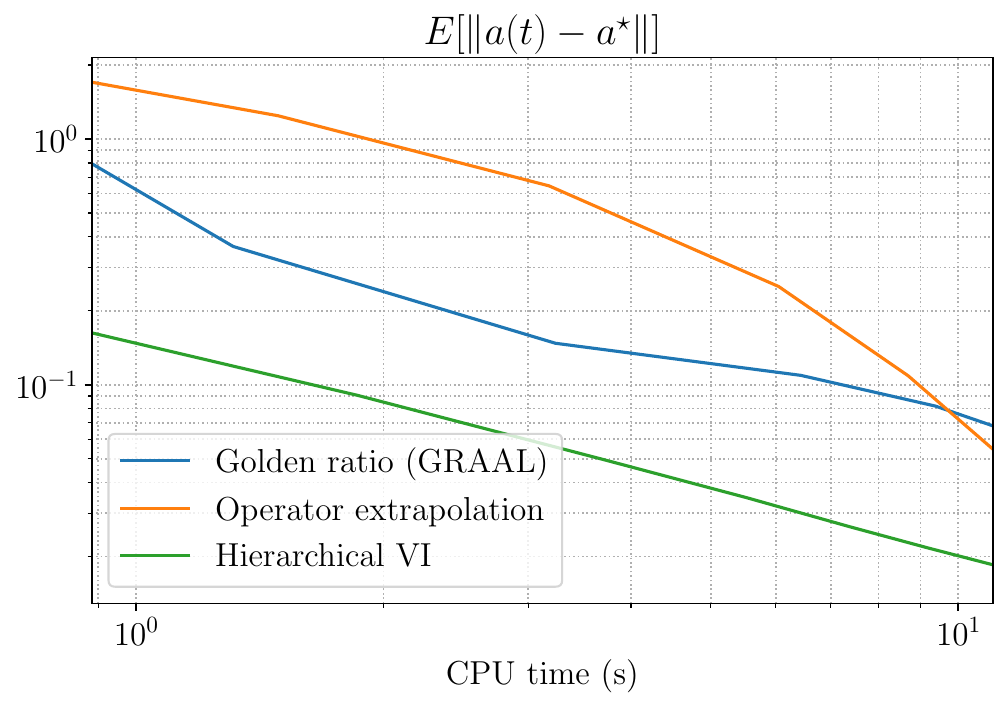}
\caption{Centralized GNE seeking, nonlinear constraints: 
Distance to the true v-GNE $\ab^\star$ as a function of CPU time.}
\label{fig:centralized_nonlinear_constraints}
\end{figure}

We test the centralized algorithm of Section~\ref{sec:centralized} using $N=10$ players with $d_i=3$, $i=1,\dots,N$ (i.e., $30$ variables in total), coupled through $500$ nonlinear shared constraints of the form~\eqref{eq:sim-coupling}, $3$ of which are active at the v-GNE. We use the Huber profile of Example~\ref{ex:huber} with $\omega=0.03$, $\varepsilon_t=0.01/(t+1)^{2/3}$, $\eta=10\gamma$, and $\gamma=2/\hat L_{\Pb}$, where $\hat L_{\Pb}$ is a numerical estimate of the Lipschitz constant of $\Pb$. This stepsize sits at the boundary of the admissible range $(0,2/L_{\Pb})$ of Theorems~\ref{th:regularized} and~\ref{th:centralized}. On this family we compare the proposed algorithm against the Golden Ratio Algorithm (GRAAL)~\cite{malitsky2020golden} on the lifted primal--dual VI~\eqref{eq:liftedVI}, and Operator Extrapolation~\cite{Lan2022} using projection steps onto the nonlinear feasible set at each iteration and applied directly to the non-lifted variational GNE problem. As shown in Fig.~\ref{fig:centralized_nonlinear_constraints}, the proposed hierarchical VI algorithm attains the smallest distance to the v-GNE. Note that, while our algorithm
requires only projections onto local box constraints, Operator Extrapolation solves an auxiliary optimization problem to project onto the nonlinear coupled feasible set at every iteration, and GRAAL operates in the lifted primal--dual space of dimension $30+500=530$.

\subsection{Distributed, nonlinear shared constraints}
\begin{figure}[t]
\centering
\includegraphics[width=0.9\linewidth]{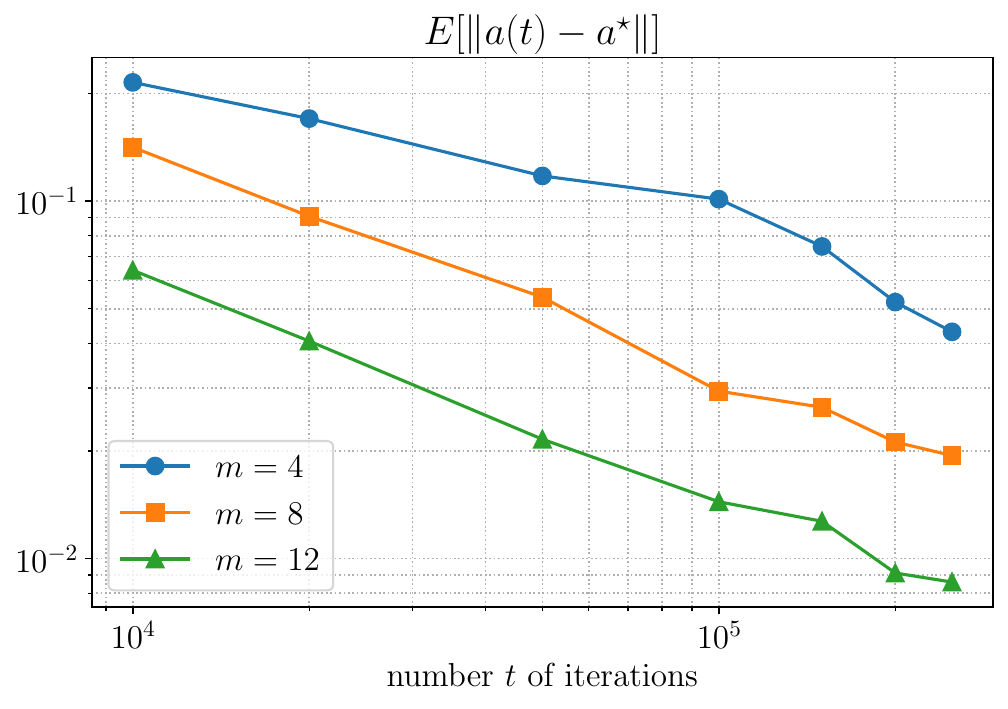}
\caption{Distributed GNE seeking, nonlinear constraints: Average distance to the true v-GNE 
for a varying number of shared nonlinear constraints.}
\label{fig:distributed_nonlinear_constraints}
\end{figure}

% example_hvi_distr_nonlinear.py
We next test the distributed hierarchical VI algorithm (Algorithm~\ref{alg:1}) on instances generated as described above with $N=10$ agents with scalar actions, i.e., $d_i=1$, $i\in[N]$, no linear shared constraints, and $m\in\{4,8,12\}$ nonlinear shared constraints split evenly between log-sum-exp and quadratic terms, $m_a=4$ of which are active at the v-GNE $\ab^\star$. %For each random seed, the three instances share the same costs, local sets, and v-GNE, and only differ in the shared constraints. 
The agents communicate over a fixed random connected graph with $22$ out of $45$ possible edges, using Metropolis--Hastings weights $w_{ij}=1/\max\{|\mathcal N_i|,|\mathcal N_j|\}$ for $j\in\mathcal N_i\setminus\{i\}$ and $w_{ii}=1-\sum_{j\neq i}w_{ij}$, which gives $\sigma\approx 0.78$ in~\eqref{eq:sigma}. We use the Huber profile of Example~\ref{ex:huber} with $\omega=0.03$, $\varepsilon_t=0.01/(t+1)^{2/3}$, $\gamma=2\cdot 0.9999/\hat L_{\Pb}$, where $\hat L_{\Pb}$ is a numerical estimate of the Lipschitz constant of $\Pb$, so that $\gamma$ lies just inside the admissible range $(0,2/L_{\Pb}^{\Qb})$ of Theorem~\ref{th:one-level-constant-step}, $\eta=10\gamma$, $\Qb^j=\Ab^j$, and $\xb_i(0)=\boldsymbol 0$ for all $i\in[N]$. Fig.~\ref{fig:distributed_nonlinear_constraints} reports the distance $\|\ab(t)-\ab^\star\|$ of the agents' actions $\ab(t)=[\xb_1^1(t),\ldots,\xb_N^N(t)]$ to the v-GNE, averaged over $20$ random instances. 

\subsection{Distributed, linear shared constraints}\label{sec:sim-distr-lin}
%example_hvi_beats_pppa_converged.py
To compare against the distributed preconditioned proximal-point algorithm (PPPA) of~\cite{bianchi2022fast} and its accelerated (over-relaxed) version (PPPA, acc), in which each agent maintains, besides its estimates of the joint action, a local copy of the dual variable $\bla\in\R^m$ that is driven to consensus over the unweighted path graph, we finally consider games with linear shared constraints $\gb(\ab)=A\ab-\bb\leq\boldsymbol 0$. The GNEPs are generated as described in Section~\ref{sec:sim-centralized-nonlinear} with $A$ having i.i.d.\ standard Gaussian entries and with large local boxes $\Ab^i$, whose bounds are obtained by subtracting from and adding to each entry of $\ab^{\star,i}$ numbers uniformly distributed in $[10,20]$. Each cost $J^i$ is linear-quadratic in the joint action $\ab$. We consider $N=5$ agents with $d_i=20$, $i\in[N]$ (i.e., $100$ variables in total), coupled through $m=300$ linear shared constraints, $100$ of which are active at the v-GNE. 

\begin{figure}[t]
\centering
\includegraphics[width=0.9\linewidth]{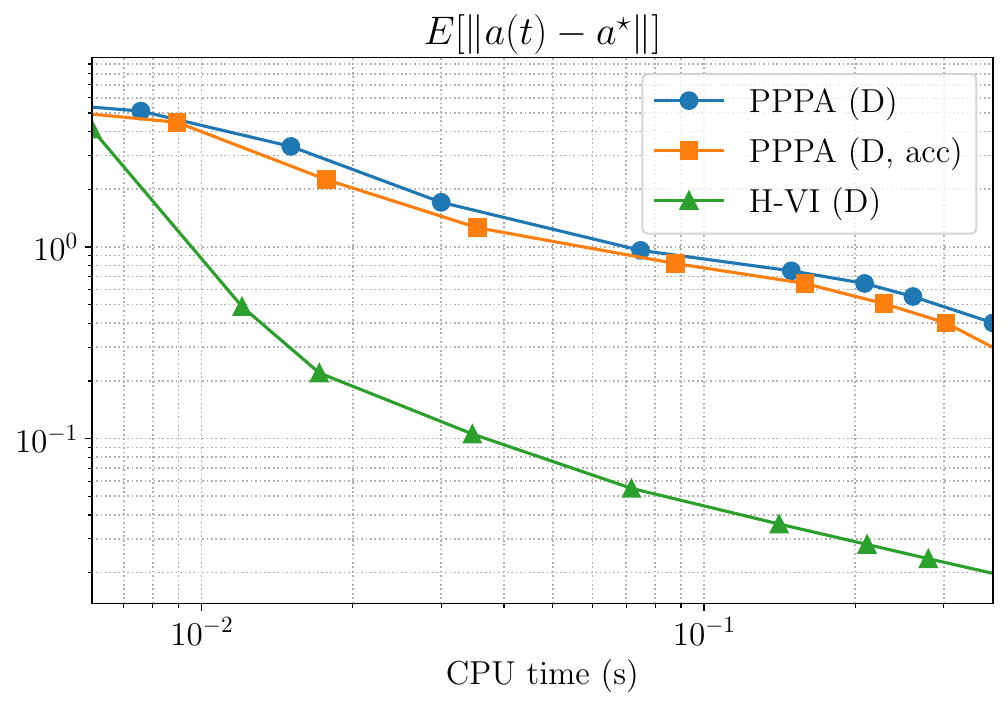}
\caption{Distributed GNE seeking, linear constraints: 
Distance to the true v-GNE as a function of CPU time.}
\label{fig:distributed_linear_constraints}
\end{figure}

The agents communicate over a path graph, the sparsest connected topology. The distributed hierarchical VI algorithm uses Metropolis--Hastings weights, which give $\sigma\approx 0.87$ in~\eqref{eq:sigma}, the Huber profile of Example~\ref{ex:huber} with $\omega=0.3$, $\varepsilon_t=1/(t+1)^{2/3}$, $\gamma=2\cdot 0.9999/\|A\|^2$ (note that $L_{\Pb}\leq\|A\|^2$), $\eta=10\gamma$, and $\xb_i(0)=\boldsymbol 0$. The PPPA stepsizes $\tau_i$, $\delta_i$, and $\nu$ are set to $0.9999$ times the upper bounds given by the sufficient conditions in~\cite[Lemma~5]{bianchi2022fast}, so that these conditions hold strictly, while $\alpha$ is set to $0.9999$ times the bound in~\cite[Lemma~3]{bianchi2022fast}, which we multiplied by the boost factor $10^6$ to overcome the large conservativeness of such a theoretical bound (hence convergence of PPPA is not formally guaranteed in this setting). As shown in Fig.~\ref{fig:distributed_linear_constraints}, which reports the average distance to the v-GNE over $20$ instances against the average execution time, the distributed hierarchical VI algorithm outperforms both versions of PPPA. This behavior is due to the fact that, in PPPA, each agent has to maintain and reach consensus on a local copy of the $300$-dimensional dual variable in addition to its estimates of the joint action, and this consensus is slow over the path graph. In contrast, the proposed algorithm exchanges only the estimates of the joint action and involves no dual variables.

\section{Conclusion}\label{sec:conclusion}
We addressed the computation of variational GNE in games with strongly monotone pseudo-gradients and nonlinear convex shared constraints. We showed that the v-GNE problem is exactly equivalent to a hierarchical VI, whose lower level is defined by the gradient of a smooth convex penalty of the constraint violation over the local action sets. This reformulation removes the need for projections onto the coupled feasible set and for dual variables. Building on it, we proposed a centralized hierarchical regularization scheme that requires only one projection onto the local action sets per iteration, and proved its last-iterate convergence to the v-GNE with rate $O(1/t^{\scriptscriptstyle 2/3})$ under a local error bound condition (implied, e.g., by Slater's condition) and $O(1/t)$ under a global error bound condition. We then extended the scheme to the partial-decision information setting and obtained a distributed algorithm with constant stepsize, in which the agents exchange only local estimates of the joint action with their neighbors, and proved its convergence to the v-GNE with rate $O(1/t)$. The numerical experiments %on instances inspired by multi-energy systems with nonlinear network constraints and on distributed problems with linear shared constraints 
showed that the proposed methods reach a higher accuracy than the projection-based and lifted primal--dual methods for the same computational time.

Future research directions include the extension of the approach to merely monotone games, to time-varying and directed communication graphs, and to settings with stochastic or payoff-based information on the agents' costs.

\appendix
\section{Appendix}\label{app:aux}
{\subsection{Properties of the Penalty Function}\label[appendix]{app:penalty}
This appendix collects the regularity properties of the penalty function $\Psi$ defined in~\eqref{eq:pen} and of its gradient $\Pb=\nabla\Psi$ in~\eqref{eq:LL}.
\begin{lem}\label{lem:P-Lip}
Let \cref{assum:convex}(ii) hold and let $\psi$ be an admissible penalty profile, with $\psi'$ and $[\cdot]_+$ acting componentwise on vectors.

(i) The function $\Psi$ is continuously differentiable and
\[
\Pb(\ab)=\nabla\gb(\ab)^\top\psi'([\gb(\ab)]_+).
\]

(ii) Let $X\subseteq\R^D$ be convex and let $\psi'$ be $L_\psi$-Lipschitz continuous. Suppose that for all $\ab,\bb\in X$
\begin{align*}
&\|\nabla\gb(\ab)\|\le G_{\gb},\qquad \|\psi'([\gb(\ab)]_+)\|\le\bar\psi,\cr
&\|\nabla\gb(\ab)-\nabla\gb(\bb)\|\le L_{\gb}\|\ab-\bb\|.
\end{align*}
Then $\Pb$ is $L_{\Pb}$-Lipschitz continuous on $X$ with $L_{\Pb}=L_\psi G_{\gb}^2+\bar\psi L_{\gb}$. If $\gb$ is affine, then $L_{\gb}=0$ and the bound $\bar\psi$ is not needed.
\end{lem}
\begin{proof}
(i) Let $\theta(s):=\psi([s]_+)$. Then $\theta'(s)=\psi'([s]_+)$ for $s\neq0$. Since $\psi'$ is nondecreasing ($\psi$ being convex), the mean value theorem yields, for $s\neq0$,
\[
\frac{|\theta(s)-\theta(0)-0\cdot s|}{|s|}
=\frac{\psi([s]_+)}{|s|}
\le\frac{\psi'([s]_+)[s]_+}{|s|}
\le\psi'([s]_+),
\]
which vanishes as $s\to0$ owing to continuity of $\psi'$ and $\psi'(0)=0$. Hence $\theta$ is continuously differentiable with $\theta'(s)=\psi'([s]_+)$, and the expression for $\Pb$ follows from the chain rule applied to $\Psi(\ab)=\sum_{j=1}^m\theta(g_j(\ab))$.

(ii) For $\ab,\bb\in X$ we have
\begin{align*}
\|\Pb(\ab)-\Pb(\bb)\|
\le{}&\|\nabla\gb(\ab)^\top\bigl(\psi'([\gb(\ab)]_+)-\psi'([\gb(\bb)]_+)\bigr)\|\cr
&+\|\bigl(\nabla\gb(\ab)-\nabla\gb(\bb)\bigr)^\top\psi'([\gb(\bb)]_+)\|\cr
\le{}&L_\psi G_{\gb}\|\gb(\ab)-\gb(\bb)\|+\bar\psi L_{\gb}\|\ab-\bb\|\cr
\le{}&(L_\psi G_{\gb}^2+\bar\psi L_{\gb})\|\ab-\bb\|,
\end{align*}
where the second inequality uses $\|[\ub]_+-[\vb]_+\|\le\|\ub-\vb\|$, and the last one follows from the mean value inequality on the segment $[\ab,\bb]\subseteq X$.
\end{proof}}

{\subsection{An Error Bound under Slater's Condition}\label[appendix]{app:slater}
\begin{lem}\label{lem:slater-eb}
Let Assumption~\ref{assum:convex}{(ii)} hold, let $\bar\ab\in\Ab$ satisfy Slater's condition, i.e., $g_j(\bar\ab)<0$ for all $j\in[m]$, and let $\bar s=\min_{j\in[m]}(-g_j(\bar\ab))>0$. Then
\[
\operatorname{dist}(\ab,\Ab\cap \Cb)\le\frac{\|\ab-\bar\ab\|}{\bar s}\,\|[\gb(\ab)]_+\|,\qquad\forall\ab\in\Ab.
\]
In particular, \cref{assum:local-eb} holds: for every bounded set $B\subseteq\R^D$ there exists $\kappa_B>0$ such that $\operatorname{dist}(\ab,\Ab\cap \Cb)\le\kappa_B\|[\gb(\ab)]_+\|$ for all $\ab\in B\cap\Ab$.
\end{lem}
\begin{proof}
Let $\ab\in\Ab$ and $r=\max_{j\in[m]}[g_j(\ab)]_+$. If $r=0$, then $\ab\in \Ab\cap \Cb$ and the claim holds. Otherwise, let $\lambda=\bar s/(\bar s+r)\in(0,1)$ and $\yb=\lambda\ab+(1-\lambda)\bar\ab$. Then $\yb\in\Ab$ by convexity of $\Ab$, and convexity of $g_j$ gives $g_j(\yb)\le\lambda r-(1-\lambda)\bar s=0$ for all $j\in[m]$, so that $\yb\in \Ab\cap \Cb$. Hence
\[
\operatorname{dist}(\ab,\Ab\cap \Cb)\le\|\ab-\yb\|=\frac{r\,\|\ab-\bar\ab\|}{\bar s+r}\le\frac{\|\ab-\bar\ab\|}{\bar s}\,r,
\]
and $r\le\|[\gb(\ab)]_+\|$. The second claim follows with $\kappa_B=\sup_{\ab\in B}\|\ab-\bar\ab\|/\bar s$.
\end{proof}}

{\subsection{A Chung-Type Lemma}\label[appendix]{app:chung}
The following lemma is a version of Chung's lemma~\cite{Chung} (see also~\cite[Ch.~2]{polyak}). We include a short proof for completeness.
\begin{lem}\label{lem:chung}
Let $\{u_t\}\subseteq[0,\infty)$ satisfy
\[
u_{t+1}\le\Bigl(1-\frac{a}{(t+1)^p}\Bigr)u_t+\frac{b}{(t+1)^q},\qquad\forall t\ge T,
\]
where $a,b>0$, $p\in(0,1]$, and $p<q\le p+1$. Then the following hold.
\begin{enumerate}
\item[(i)] If $p<1$, then $u_t=O(t^{-(q-p)})$.
\item[(ii)] If $p=1$ and $a>q-1$, then $u_t=O(t^{-(q-1)})$.
\item[(iii)] If $p=1$, then $u_t\to0$.
\end{enumerate}
\end{lem}
\begin{proof}
Let $x_t=1/(t+1)\in(0,1]$, and let $s=q-p$ in case (i), $s=q-1$ in case (ii), and $s=\min\{a,q-1\}/2$ in case (iii). In each case $0<s\le q-p\le1$, and $s<a$ if $p=1$. We show that $v_t=(t+1)^su_t$ is bounded, which gives $u_t=O(t^{-s})$ and proves all three claims.

Since $(t+2)/(t+1)=1+x_t$, multiplying the recursion by $(t+2)^s$ gives
\[
v_{t+1}\le(1+x_t)^s(1-ax_t^p)\,v_t+b\,(1+x_t)^s\,x_t^{q-s}.
\]
Since $s\le1$, we have $(1+x_t)^s\le1+sx_t\le2$, and since $q-s\ge p$, we have $x_t^{q-s}\le x_t^p$. Hence, for all $t$ large enough that $ax_t^p\le1$,
\[
v_{t+1}\le(1+sx_t-ax_t^p)\,v_t+2b\,x_t^p.
\]
Let $a'=a/2$ if $p<1$ and $a'=a-s$ if $p=1$. Then $sx_t-ax_t^p\le-a'x_t^p$ for all large $t$. For $p<1$ this holds because $x_t^{1-p}\to0$, and for $p=1$ it holds with equality. With $\theta_t=a'x_t^p\in(0,1]$ and $K=2b/a'$, we thus have, for all $t\ge T'$ and some $T'\ge T$,
\[
v_{t+1}\le(1-\theta_t)\,v_t+\theta_tK.
\]
The right-hand side is a convex combination of $v_t$ and $K$, so $v_{t+1}\le\max\{v_t,K\}$. By induction, $v_t\le\max\{v_{T'},K\}$ for all $t\ge T'$.
\end{proof}}

% Proofs of the auxiliary lemmas of Section~\ref{sec:conv_distr}
% (kept as subsections of the single appendix; ieeecolor prints every appendix \section as an untitled "Appendix" heading)

\subsection{Proof of Lemma~\ref{lem:one-level-owned-disagreement}}
Let
$\ab=\ab(\xb)$,
$\overline{\xb} =
\frac{1}{N}\sum_{i=1}^{N}\xb_i$.
For every block $k\in[N]$, define
\[
\overline{\xb}^{\,k}
=
\frac{1}{N}\sum_{i=1}^{N}\xb_i^k.
\]
Since
$\ab^k=\xb_k^k$,
we obtain
\begin{align*}
\sum_{i=1}^{N}\|\xb_i^k-\ab^k\|^2
&=
\sum_{i=1}^{N}\|\xb_i^k-\xb_k^k\|^2
\cr
&=
\sum_{i=1}^{N}
\left\|
(\xb_i^k-\overline{\xb}^{\,k})
-
(\xb_k^k-\overline{\xb}^{\,k})
\right\|^2
\cr
&=
\sum_{i=1}^{N}
\|\xb_i^k-\overline{\xb}^{\,k}\|^2
+
N\|\xb_k^k-\overline{\xb}^{\,k}\|^2
\cr
&\leq
(N+1)
\sum_{i=1}^{N}
\|\xb_i^k-\overline{\xb}^{\,k}\|^2.
\end{align*}
Summing over $k\in[N]$ proves
\[
\left\|
\xb
-
\boldsymbol 1_N\otimes \ab
\right\|^2
\leq
(N+1)
\left\|
\xb
-
\boldsymbol 1_N\otimes\overline{\xb}
\right\|^2.
\]
Since
$\xb
-
\boldsymbol 1_N\otimes\overline{\xb}
=
\xb^\perp$,
the claim follows.

\subsection{An Auxiliary Gradient-Dominance Bound}
The following lemma is used only internally, in the proof of Lemma~\ref{lem:one-level-residual-control} below.
\begin{lem}
\label{lem:one-level-gradient-dominance}
Let $\psi$ be an admissible penalty profile with constants $S_0,c_\psi>0$ as in Definition~\ref{def:profile}. Under Assumptions~\ref{assum:convex}{(ii)},~\ref{assum:local-eb}, and~\ref{assum:box}, {let $\kappa_{\Ab}$ be the constant of \cref{assum:local-eb} for $B=\Ab$. Then} there exists $c'>0$, depending only on $\psi$ and on $\Ab$, such that
$
\Psi(\ab)
\le
\kappa_{\Ab}^2(c')^2\|\Pb(\ab)\|^2$,
$
\forall \ab\in\Ab.
$
Consequently,
$
\operatorname{dist}(\ab,S_L)^2
\leq
\kappa_{\Ab}^4(c')^4\|\Pb(\ab)\|^2$,
$\forall \ab\in\Ab$.
\end{lem}
\begin{proof}
Since $\Ab$ is bounded (Assumption~\ref{assum:box}), Assumption~\ref{assum:local-eb} with $B=\Ab$ gives $\operatorname{dist}(\ab,S_L)\le\kappa_{\Ab}\|[\gb(\ab)]_+\|$ for every $\ab\in\Ab$. %we use this $\kappa_{\Ab}$ below.
Since $\Ab$ is closed (Assumption~\ref{assum:convex}{(ii)}) and bounded (Assumption~\ref{assum:box}), it is compact, so $\Psi_{\max}=\max_{\ab\in\Ab}\Psi(\ab)<\infty$ is attained ($\Psi$ being continuous). Set
\[
c'=\max\left(\frac{1}{\sqrt{c_\psi}},\ \frac{\sqrt{\Psi_{\max}}}{c_\psi S_0}\right).
\]
Let $\ab\in\Ab$ and write $r=\|[\gb(\ab)]_+\|$. By Definition~\ref{def:profile} applied coordinatewise and Lemma~\ref{lem:norm-hybrid}, recalling~\eqref{eq:pen}, $\Psi(\ab)\ge c_\psi\, r\min(r,S_0)$; since $s\mapsto s\min(s,S_0)$ is increasing, this inverts to
\[
r \le \max\!\left(\sqrt{\Psi(\ab)/c_\psi},\ \frac{\Psi(\ab)}{c_\psi S_0}\right).
\]
As $\Psi(\ab)\le\Psi_{\max}$, the second term satisfies $\frac{\Psi(\ab)}{c_\psi S_0}=\frac{\sqrt{\Psi(\ab)}}{c_\psi S_0}\sqrt{\Psi(\ab)}\le\frac{\sqrt{\Psi_{\max}}}{c_\psi S_0}\sqrt{\Psi(\ab)}$, so that, combined with the first term $\sqrt{\Psi(\ab)/c_\psi}=\frac{1}{\sqrt{c_\psi}}\sqrt{\Psi(\ab)}$,
\begin{align}\label{eq:L1_0}
\|[\gb(\ab)]_+\| = r \le c'\sqrt{\Psi(\ab)}.
\end{align}

Let $\sb=\Pi_{S_L}(\ab)$. Since \(\Psi(\sb)=0\), convexity of $\Psi$ implies
\[
0
=
\Psi(\sb)
\geq
\Psi(\ab)+\langle \Pb(\ab),\sb-\ab\rangle,
\]
so that
\[
\Psi(\ab)
\le
\|\Pb(\ab)\|
\operatorname{dist}(\ab,S_L).
\]
By the error bound established above and~\eqref{eq:L1_0},
\begin{align}\label{eq:L1_1}
\operatorname{dist}(\ab,S_L)
\le
\kappa_{\Ab}\|[\gb(\ab)]_+\|
\le
\kappa_{\Ab} c'\sqrt{\Psi(\ab)}.
\end{align}
Combining the relations above, we obtain
$
\Psi(\ab)
\le
\kappa_{\Ab} c'
\|\Pb(\ab)\|
\sqrt{\Psi(\ab)}$.
If $\Psi(\ab)=0$, the claims hold. Otherwise, division by
$\sqrt{\Psi(\ab)}$ and squaring yield
\[
\Psi(\ab)
\le
\kappa_{\Ab}^2(c')^2\|\Pb(\ab)\|^2.
\]
Finally, using~\eqref{eq:L1_1}, we conclude
\[
\operatorname{dist}(\ab,S_L)^2
\leq
\kappa_{\Ab}^2(c')^2\Psi(\ab)
\leq
\kappa_{\Ab}^4(c')^4\|\Pb(\ab)\|^2.
\]
\end{proof}

\subsection{Proof of Lemma~\ref{lem:one-level-residual-control}}
Set $\ab=\ab(\xb)\in\Ab\subseteq \Qb$.
By Lemma~\ref{lem:one-level-owned-disagreement},
\[
\left\|
\xb
-
\boldsymbol 1_N\otimes \ab
\right\|^2
\leq
(N+1)\|\xb^\perp\|^2.
\]
Since $W\boldsymbol 1_N=\boldsymbol 1_N$ (Assumption~\ref{assum:comm}), we have
$\Wb(\boldsymbol 1_N\otimes \ab)=(W\otimes I_D)(\boldsymbol 1_N\otimes \ab)=(W\boldsymbol 1_N)\otimes \ab=\boldsymbol 1_N\otimes \ab$, so that
\[
\zb
-
\boldsymbol 1_N\otimes \ab
=
\Wb\xb-\Wb(\boldsymbol 1_N\otimes \ab)
=
\Wb
\left(
\xb
-
\boldsymbol 1_N\otimes \ab
\right).
\]
Moreover, $W$, being symmetric and doubly stochastic, has eigenvalues in $[-1,1]$, so $\|\Wb\|=\|W\otimes I_D\|=\|W\|\leq1$. Combining this inequality with the previous equality, we obtain 
\begin{align}\label{eq:zDisagr}
\left\|
\zb
-
\boldsymbol 1_N\otimes \ab
\right\|^2
\leq
\|\Wb\|^2\left\|\xb-\boldsymbol 1_N\otimes \ab\right\|^2
\leq
(N+1)\|\xb^\perp\|^2.
\end{align}

As $\xb\in\boldsymbol{\Omega}$, we have $\zb_i\in \Qb$ for every $i\in[N]$. Since also $\ab\in \Qb$, Assumption~\ref{assum:box} implies, for every $i\in[N]$,
\[
\|\Pb(\ab)\|
\leq
\|\Pb(\zb_i)\|+L_{\Pb}^{\Qb}\|\zb_i-\ab\|.
\]
Squaring via $(u+v)^2\leq2u^2+2v^2$, summing over $i\in[N]$, and dividing by $N$ gives
\begin{align*}
\|\Pb(\ab)\|^2
&\leq
\frac{2}{N}\sum_{i=1}^N\|\Pb(\zb_i)\|^2
+
\frac{2(L_{\Pb}^{\Qb})^2}{N}\sum_{i=1}^N\|\zb_i-\ab\|^2\cr
&=
\frac{2}{N}
\|\Pb(\zb)\|^2
+
\frac{2(L_{\Pb}^{\Qb})^2}{N}
\left\|
\zb
-
\boldsymbol 1_N\otimes \ab
\right\|^2,
\end{align*}
where we used $\Pb(\zb)=[\Pb(\zb_1),\ldots,\Pb(\zb_N)]$ and $\|\zb-\boldsymbol 1_N\otimes \ab\|^2=\sum_{i=1}^N\|\zb_i-\ab\|^2$.
Substituting $\|\zb-\boldsymbol 1_N\otimes \ab\|^2\leq(N+1)\|\xb^\perp\|^2$ from~\eqref{eq:zDisagr} and using $\|\xb^\perp\|^2\leq\mathcal{R}(\xb)$, $\|\Pb(\zb)\|^2\leq\mathcal{R}(\xb)$ (as $\mathcal{R}(\xb)=\|\xb^\perp\|^2+\|\Pb(\zb)\|^2$), we obtain
\begin{align*}
\|\Pb(\ab)\|^2
&\leq
\frac{2}{N}\|\Pb(\zb)\|^2
+
\frac{2(L_{\Pb}^{\Qb})^2(N+1)}{N}\|\xb^\perp\|^2\leq
C_{\Pb}\mathcal{R}(\xb),
\end{align*}
where $C_{\Pb}=2\max(1,(L_{\Pb}^{\Qb})^2(N+1)/N)$.
Thus, Lemma~\ref{lem:one-level-gradient-dominance} together with the inequality above
 gives
\begin{align}\label{eq:distS}
\operatorname{dist}(\ab,S_L)^2
\leq
\kappa_{\Ab}^4(c')^4C_{\Pb}\mathcal{R}(\xb).
\end{align}
Next, for every $\sb\in S_L$,
\begin{align*}
\|\xb-
\boldsymbol 1_N\otimes \sb
\|
&\leq
\left\|
\xb
-
\boldsymbol 1_N\otimes \ab
\right\|
+
\left\|\boldsymbol 1_N\otimes \ab-\boldsymbol 1_N\otimes \sb\right\|\cr
&=
\left\|
\xb
-
\boldsymbol 1_N\otimes \ab
\right\|
+
\sqrt{N}\|\ab-\sb\|,
\end{align*}
where we used $\|\boldsymbol 1_N\otimes \vb\|^2=N\|\vb\|^2$ for $\vb=\ab-\sb$. Since $\boldsymbol{S}_L=\{\boldsymbol 1_N\otimes\sb:\sb\in S_L\}$, taking the infimum over $\sb\in S_L$ on both sides gives
\[
\operatorname{dist}(\xb,\boldsymbol{S}_L)
\leq
\left\|
\xb
-
\boldsymbol 1_N\otimes \ab
\right\|
+
\sqrt{N}\operatorname{dist}(\ab,S_L).
\]
Squaring via $(u+v)^2\leq2u^2+2v^2$, we obtain
\[
\begin{aligned}
\operatorname{dist}
\left(
\xb,
\boldsymbol{S}_L
\right)^2
&\leq
2\left\|\xb-\boldsymbol 1_N\otimes \ab\right\|^2
+
2N\operatorname{dist}(\ab,S_L)^2
\\
&\leq
2(N+1)\|\xb^\perp\|^2
+
2N\kappa_{\Ab}^4(c')^4C_{\Pb}\mathcal{R}(\xb)
\\
&\leq
\Bigl[2(N+1)+2N\kappa_{\Ab}^4(c')^4C_{\Pb}\Bigr]\mathcal{R}(\xb),
\end{aligned}
\]
where in the first inequality we used Lemma~\ref{lem:one-level-owned-disagreement} and~\eqref{eq:distS}, and in the last inequality we used $\|\xb^\perp\|^2\leq\mathcal{R}(\xb)$.
Thus,
\[
\operatorname{dist}
\left(
\xb,
\boldsymbol{S}_L
\right)^2
\leq
C_R\mathcal{R}(\xb),
\]
where $C_R=2(N+1)+2N\kappa_{\Ab}^4(c')^4C_{\Pb}$.
Finally,
\[
\begin{aligned}
\|
\xb
-
\boldsymbol 1_N\otimes \sb
\|^2
&\leq
2
\left\|
\xb
-
\boldsymbol 1_N\otimes \ab
\right\|^2
+
2N\|\ab-\sb\|^2
\\
&\leq
2(N+1)\|\xb^\perp\|^2
+
2N\|\ab-\sb\|^2.
\end{aligned}
\]
Thus one may take
$K_R=2(N+1)$,
and $K_a=2N$.

\bibliographystyle{plain}
\bibliography{srtrMonGames_ref}

\vfill

\end{document}